\documentclass[a4paper,10pt]{amsproc}
\usepackage{amsfonts, amsmath, amssymb, graphicx, mathrsfs, tikz, xcolor, hyperref}
\usepackage{chngcntr}
\usepackage[a4paper, total={6in, 10in}]{geometry}
\hypersetup{
	colorlinks   = true, 
	urlcolor     = blue, 
	linkcolor    = green, 
	citecolor   = red 
}
\tikzset{Bullet/.style={fill=black,draw,color=#1,circle,minimum size=2pt,scale=0.45}}
\usetikzlibrary{shapes.geometric}
\usetikzlibrary{positioning, tikzmark}
\usepackage[all]{xy}
\theoremstyle{plain}

\newenvironment{section*}[2][E]{
	\section*{#2}
	\renewcommand\thesection{\thechapter.#1}
	\setcounter{theorem}{0}}{}
\newtheorem{theorem}{Theorem}[section]
\newtheorem{lemma}[theorem]{Lemma}

\theoremstyle{definition}
\newtheorem{definition}[theorem]{Definition}

\newtheorem{remark}[theorem]{Remark}
\newtheorem{question}[theorem]{Question}
\newtheorem{counter example}[theorem]{Counter Example}
\newtheorem{notation}[theorem]{Notation}
\newtheorem{corollary}[theorem]{Corollary}
\newtheorem{example}[theorem]{Example}

\numberwithin{equation}{section}

\begin{document}
	\Large{\title[$U$-topology and $m$-topology on the ring of Measurable Functions]{$U$-topology and $m$-topology on the ring of Measurable Functions, generalized and revisited}
		\author[P. Nandi]{Pratip Nandi}
		\address{Department of Pure Mathematics, University of Calcutta, 35, Ballygunge Circular Road, Kolkata 700019, West Bengal, India}
		\email{pratipnandi10@gmail.com}
		\author[A. Deb Ray]{Atasi Deb Ray}
		\address{Department of Pure Mathematics, University of Calcutta, 35, Ballygunge Circular Road, Kolkata 700019, West Bengal, India}
		\email{debrayatasi@gmail.com}
		\author[S.K. Acharyya]{Sudip Kumar Acharyya}
		\address{Department of Pure Mathematics, University of Calcutta, 35, Ballygunge Circular Road, Kolkata 700019, West Bengal, India}
		\email{sdpacharyya@gmail.com}
		
		\thanks{The first author thanks the CSIR, New Delhi – 110001, India, for financial support.}
		
		\begin{abstract}
			Let $\mathcal{M}(X,\mathcal{A})$ be the ring of all real valued measurable functions defined over the measurable space $(X,\mathcal{A})$. Given an ideal $I$ in $\mathcal{M}(X,\mathcal{A})$ and a measure $\mu:\mathcal{A}\to[0,\infty]$, we introduce the $U_\mu^I$-topology and the $m_\mu^I$-topology on $\mathcal{M}(X,\mathcal{A})$ as generalized versions of the topology of uniform convergence or the $U$-topology and the $m$-topology on $\mathcal{M}(X,\mathcal{A})$ respectively. With $I=\mathcal{M}(X,\mathcal{A})$, these two topologies reduce to the $U_\mu$-topology and the $m_\mu$-topology on $\mathcal{M}(X,\mathcal{A})$ respectively, already considered in \cite{Acharyya3}. It is realized that if $I$ is a countably generated ideal in $\mathcal{M}(X,\mathcal{A})$, then the  $U_\mu^I$-topology and the $m_\mu^I$-topology coincide if and only if $X\setminus \bigcap Z[I]$ is a $\mu$-bounded subset of $X$. Here $Z[I]$ is the filter of measurable sets on $X$ corresponding to the ideal $I$. A subset $Y$ of $X$ is called $\mu$-bounded if each function in $\mathcal{M}(X,\mathcal{A})$ is bounded on $Y$, excepting on a set of $\mu$-measure $0$. The components of $0$ in $\mathcal{M}(X,\mathcal{A})$ in the $U_\mu^I$-topology and the $m_\mu^I$-topology are realized as $I\cap L^\infty(X,\mathcal{A},\mu)$ and $I\cap L_\psi(X,\mathcal{A},\mu)$ respectively. Here $L^\infty(X,\mathcal{A},\mu)$ is the set of all functions in $\mathcal{M}(X,\mathcal{A})$ which are essentially $\mu$-bounded over $X$ and $L_\psi(X,\mathcal{A},\mu)=\{f\in \mathcal{M}(X,\mathcal{A}): ~\forall g\in\mathcal{M}(X,\mathcal{A}), f.g\in L^\infty(X,\mathcal{A},\mu)\}$. It is established that an ideal $I$ in $\mathcal{M}(X,\mathcal{A})$ is dense in the $U_\mu$-topology if and only if it is dense in the $m_\mu$-topology and this happens when and only when there exists $Z\in Z[I]$ such that $\mu(Z)=0$. Furthermore, it is proved that $I$ is closed in $\mathcal{M}(X,\mathcal{A})$ in the $m_\mu$-topology if and only if it is a $Z_\mu$-ideal in the sense that if $f\equiv g$ almost everywhere on $X$ with $f\in I$ and $g\in\mathcal{M}(X,\mathcal{A})$, then $g\in I$. Finally two more topologies on $\mathcal{M}(X,\mathcal{A})$ viz. the $U_{\mu,F}^I$-topology and the $m_{\mu,F}^I$-topology, coarser than the $U_\mu^I$-topology and the $m_\mu^I$-topology respectively are introduced and a few relevant properties are investigated thereon.
		\end{abstract}
	
	\maketitle
	\section{Introduction}
	Our starting point is a triple $(X,\mathcal{A},\mu)$, here $X$ is a non-empty set, $\mathcal{A}$, a $\sigma$-algebra of subsets of $X$ and $\mu:\mathcal{A}\to [0,\infty]$, a measure. This triple is familiarly called a measure space. A function $g:X\to\mathbb{R}$ is called a measurable function in this measure space if for any open set $U$ in $\mathbb{R}$, $g^{-1}(U)$ is a member of $\mathcal{A}$. The members of $\mathcal{A}$ are often called measurable sets. It is well-known that the aggregate $\mathcal{M}(X,\mathcal{A})$ of all (real valued) measurable functions makes a commutative lattice ordered ring with unity, provided the relevant compositions are defined pointwise on $X$. A number of authors have already worked on a few problems and their possible ramifications connected with the ring $\mathcal{M}(X,\mathcal{A})$. The reader may consult the articles~\cite{Acharyya1}, \cite{Acharyya2}, \cite{Acharyya3} in this regard. Let for $f\in\mathcal{M}(X,\mathcal{A})$, $\epsilon>0$, $u$, a positive unit in the ring $\mathcal{M}(X,\mathcal{A})$, $U_\mu(f,\epsilon)=\{g\in\mathcal{M}(X,\mathcal{A}):\text{there exists a set }E_g\in\mathcal{A}\text{ with }\mu(E_g)=0\text{ such that }\sup\limits_{x\in X\setminus E_g}|f(x)-g(x)|<\epsilon\}$ and $m_\mu(f,u)=\{g\in\mathcal{M}(X,\mathcal{A}):\text{there exists a }\mu\text{-null set }F\text{ with }|f(x)-g(x)|<u(x)\text{ for each }x\in X\setminus F\}$ [$F$ is called $\mu$-null set if $F\in\mathcal{A}$ and $\mu(F)=0$]. Then the families: $\{U_\mu(f,\epsilon):f\in\mathcal{M}(X,\mathcal{A}), \epsilon>0\}$ and $\{m_\mu(f,u):f\in\mathcal{M}(X,\mathcal{A}), u,\text{ a positive unit in }\mathcal{M}(X,\mathcal{A})\}$ are open bases respectively for two topologies viz., $U_\mu$-topology and $m_\mu$-topology on the ring $\mathcal{M}(X,\mathcal{A})$. A number of problems related to these two topologies are already addressed in the article~\cite{Acharyya3}. Incidentally a generalized version of each of these two topologies via an ideal $I$ in $\mathcal{M}(X,\mathcal{A})$ is initiated in the communicated paper~\cite{Acharyya4}, where a few problems relevant thereon are studied. In the present paper, we offer yet another generalization of $U_\mu$-topology and $m_\mu$-topology on $\mathcal{M}(X,\mathcal{A})$, this time by a different technique. To be precise let $I$ be an ideal in $\mathcal{M}(X,\mathcal{A})$, not necessarily a proper ideal. For the above choices of $f, \epsilon, u$, let $U_\mu(f,I,\epsilon)=\{g\in\mathcal{M}(X,\mathcal{A}):g-f\in I\text{ and } \sup\limits_{x\in X\setminus E_g}|f(x)-g(x)|<\epsilon\text{ for some }\mu\text{-null set }E_g\}$ and $m_\mu(f,I,u)=\{g\in\mathcal{M}(X,\mathcal{A}):g-f\in I\text{ and }|f(x)-g(x)|<u(x)\text{ for all }x\text{, almost everywhere with respect to the measure }\mu\text{ on }X\}$. It is easy to verify by routine computations and taking care of the finite subadditivity of the measure function $\mu:\mathcal{A}\to[0,\infty]$ that the families $\mathscr{B}_U=\{U_\mu(f,I,\epsilon):f\in\mathcal{M}(X,\mathcal{A}),\epsilon>0\}$ and $\mathscr{B}_m=\{m_\mu(f,I,u):f\in\mathcal{M}(X,\mathcal{A}), u,\text{ a positive unit in }\\\mathcal{M}(X,\mathcal{A})\}$ are open bases for two uniquely defined topologies viz., $U_\mu^I$-topology and $m_\mu^I$-topology respectively on $\mathcal{M}(X,\mathcal{A})$. With $I=\mathcal{M}(X,\mathcal{A})$, these two topologies reduce to $U_\mu$-topology and $m_\mu$-topology respectively on $\mathcal{M}(X,\mathcal{A})$ considered in~\cite{Acharyya3}. We wish to designate the $U_\mu^I$-topology (respectively $m_\mu^I$-topology) by the same notation $U_\mu^I$ (respectively $m_\mu^I$). It is not at all hard to realise on using finite subadditivity of $\mu:\mathcal{A}\to[0,\infty]$ that the topological space $(\mathcal{M}(X,\mathcal{A}), U_\mu^I)$ is an additive abelian topological group and the topological space $(\mathcal{M}(X,\mathcal{A}), m_\mu^I)$ is a topological ring. In case $(\mathcal{M}(X,\mathcal{A}), U_\mu^I)$ becomes a topological ring it is proved that $I\subseteq L^\infty(X,\mathcal{A},\mu)$, the set of all essentially bounded measurable functions on $X$. It is not known to us whether the condition $I\subseteq L^\infty(X,\mathcal{A},\mu)$ is necessary to render $(\mathcal{M}(X,\mathcal{A}), U_\mu^I)$, a topological ring. It is also realized that $(\mathcal{M}(X,\mathcal{A}), U_\mu^I)$  is a topological vector space if and only if $(\mathcal{M}(X,\mathcal{A}), m_\mu^I)$ is a topological vector space. We check that different ideals $I,J$ in $\mathcal{M}(X,\mathcal{A})$ give rise to different topologies $U_\mu ^I$ and $U_\mu ^J$ (respectively distinct $m_\mu^I$ and $m_\mu^J$). The component of $0$ in $\mathcal{M}(X,\mathcal{A})$ in the $U_\mu^I$-topology is found to be $I\cap L^\infty(X,\mathcal{A},\mu)$. The component of $0$ in $\mathcal{M}(X,\mathcal{A})$ in the $m_\mu^I$-topology is determined to be the set $I\cap L_\psi(X,\mathcal{A},\mu)$, where $L_\psi(X,\mathcal{A},\mu)=\{f\in\mathcal{M}(X,\mathcal{A}):\forall g\in\mathcal{M}(X,\mathcal{A}), f.g\in L^\infty(X,\mathcal{A},\mu)\}$. It turns out that $L_\psi(X,\mathcal{A},\mu)$ is the largest ideal in $\mathcal{M}(X,\mathcal{A})$ containing $L^\infty(X,\mathcal{A},\mu)$. We call a subset $Y$ of $(X,\mathcal{A},\mu)$ $\mu$-bounded if each function in $\mathcal{M}(X,\mathcal{A})$ is bounded on $Y$ almost everywhere. $\mu$-boundedness of certain subsets of $(X,\mathcal{A},\mu)$ is found to be equivalent to the coincidence of the  two topologies viz., $U_\mu^I$-topology and $m_\mu^I$-topology. Indeed it is established that if the ideal $I$ is countably generated, then  $U_\mu^I$-topology $=m_\mu^I$-topology if and only if $X\setminus \bigcap Z[I]$ is a $\mu$-bounded subset of $X$. A special case of this result on choosing $I=\mathcal{M}(X,\mathcal{A})$ reads: $U_\mu$-topology $=m_\mu$-topology if and only if each function in $\mathcal{M}(X,\mathcal{A})$ is essentially bounded. The following two special subrings of $\mathcal{M}(X,\mathcal{A})$ are already introduced in~\cite{Acharyya2}: $\mathcal{M}_F(X,\mathcal{A})=\{f\in\mathcal{M}(X,\mathcal{A}):\mu(X\setminus Z(f))<\infty\}$ and $\mathcal{M}_\infty(X,\mathcal{A})=\{f\in\mathcal{M}(X,\mathcal{A}):\forall \epsilon>0,\mu\{x\in X:|f(x)|\geq\epsilon\}<\infty\}=\{f\in\mathcal{M}(X,\mathcal{A}):\forall n\in\mathbb{N},\mu\{x\in X:|f(x)|\geq\frac{1}{n}\}<\infty\}$. We realize in the present article that $\mathcal{M}_\infty(X,\mathcal{A})$ is the closure of $\mathcal{M}_F(X,\mathcal{A})$ in the $U_\mu$-topology. We also determine the closure of $\mathcal{M}_F(X,\mathcal{A})$ in the $m_\mu$-topology as the set $\{f\in\mathcal{M}(X,\mathcal{A}):\forall g\in\mathcal{M}(X,\mathcal{A}), f.g\in\mathcal{M}_\infty(X,\mathcal{A})\}$. At this point we like to mention that a similar kind of topology on the ring $C(X)$ of all real valued continuous functions over a completely regular Hausdorff topological space $X$, generalizing the well known $m$-topology on $C(X)$ via an ideal $I$ in $C(X)$ is introduced in~\cite{Azarpanah5}. This topology is denoted as $m^I$-topology in~\cite{Azarpanah5}. It is proved in the same article~[\cite{Azarpanah5}, Proposition 2.2(e)] that each maximal ideal in $C(X)$ in the $m^I$-topology is closed. But in our present paper we see that a maximal ideal in $\mathcal{M}(X,\mathcal{A})$ in the $m_\mu^I$-topology need not to be closed [see the counterexample \ref{Ex5.4}]. However a certain kind of ideals viz., the $Z_\mu$-ideals in $\mathcal{M}(X,\mathcal{A})$ happens to be closed in the $m_\mu^I$-topology with the specific choice $I=\mathcal{M}(X,\mathcal{A})$. An ideal $I$ in $\mathcal{M}(X,\mathcal{A})$ is called a $Z_\mu$-ideal if whenever an $f\in\mathcal{M}(X,\mathcal{A})$ becomes zero almost everywhere on some $Z\in Z[I]$, then $f\in I$. It is realized that $Z_\mu$-ideals in $\mathcal{M}(X,\mathcal{A})$ are precisely the $m_\mu$-closed in this ring (Corollary \ref{Cor5.16}). It is observed that a maximal ideal $M$ in $\mathcal{M}(X,\mathcal{A})$ need not be a $Z_\mu$-ideal. However if $M$ contains a $Z_\mu$-ideal, then $M$ is a $Z_\mu$-ideal (Corollary \ref{Th5.199}). A set $A\in\mathcal{A}$ is called $\mu$-finite if $\mu(A)<\infty$. For an ideal $I$ in $\mathcal{M}(X,\mathcal{A})$, $f\in\mathcal{M}(X,\mathcal{A})$, $\epsilon>0$ and $u$, a positive unit in $\mathcal{M}(X,\mathcal{A})$, let $U_{\mu,F}(f,I,\epsilon)=\{g\in\mathcal{M}(X,\mathcal{A}):f-g\in I\text{ and there exists a }\mu\text{-finite set }E_g\text{ such that }\sup\limits_{x\in X\setminus E_g}|f(x)-g(x)|<\epsilon\}$ and $m_{\mu,F}(f,I,u)=\{g\in\mathcal{M}(X,\mathcal{A}):f-g\in I\text{ and there exists a }\mu\text{-finite set }E_g\text{ such that }|f(x)-g(x)|<u(x)\text{ on }X\setminus E_g\}$. It is easy to verify that the families $\mathscr{B}^F_U=\{U_{\mu,F}(f,I,\epsilon):f\in\mathcal{M}(X,\mathcal{A}),~\epsilon>0\}$ and $\mathscr{B}^F_m=\{m_{\mu,F}(f,I,u):f\in \mathcal{M}(X,\mathcal{A}),~u,\text{ a positive unit in }\mathcal{M}(X,\mathcal{A})\}$ form open bases for some topologies on $\mathcal{M}(X,\mathcal{A})$, we call them $U^I_{\mu,F}$-topology and $m^I_{\mu,F}$-topology respectively. For $I=\mathcal{M}(X,\mathcal{A})$, we simply denote these two topologies as $U_{\mu,F}$-topology and $m_{\mu,F}$-topology respectively. After a routine check it can be proved that $(\mathcal{M}(X,\mathcal{A}), U^I_{\mu,F})$ is a topological group and $(\mathcal{M}(X,\mathcal{A}), m^I_{\mu,F})$ is a topological ring for every ideal $I$ in $\mathcal{M}(X,\mathcal{A})$. Now if $X$ is $\mu$-finite, then both of these topologies reduce to the indiscrete topology on $\mathcal{M}(X,\mathcal{A})$. So for the Section \ref{Sec6}, we consider $\mu(X)=\infty$. Likewise the Section \ref{Sec2}, we state some basic facts about the $U^I_{\mu,F}$-topology and $m^I_{\mu,F}$-topology on $\mathcal{M}(X,\mathcal{A})$. Then we compute the components of $0$ in those topologies [Theorem \ref{Th6.14}(\ref{Th6.14.3}), Theorem \ref{Th6.18}]. Later we introduce a $Z_\mu$-like ideal in the finite setting viz., $Z_{\mu,F}$-ideal. By introducing the $Z_{\mu,F}$-ideal, we have proved that the class of all $Z_{\mu,F}$-ideals in $\mathcal{M}(X,\mathcal{A})$ is precisely the set of all closed ideals in the $m_{\mu,F}$-topology [Corollary \ref{Cor6.29}]. Furthermore it is proved in Theorem \ref{Th6.31}(\ref{Th6.31.3}) that a maximal ideal $M$ is closed in the $m_{\mu,F}$-topology if and only if $\mu(Z(f))=\infty$ for all $f\in M$. In section \ref{Sec7}, we illustrate some of the ideas and results related to the above mentioned four topologies by virtue of a few chosen measure spaces. Eventually we get that if $(X,\mathcal{A},\mu)$ contains an atom, then every measurable function on this measure space is essentially bounded on that atom [Lemma \ref{Lem7.4}]. An atom of the measure space is a member $A\in\mathcal{A}$ such that $\mu(A)>0$ and $A$ can not be decomposed in the manner $A=B\sqcup C$ with $B,C\in\mathcal{A}$ and $\mu(B),\mu(C)>0$. Indeed a measure space $(X,\mathcal{A},\mu)$ is called atomic if each measurable set with positive measure contains an atom. On the contrary, $(X,\mathcal{A},\mu)$ is called non-atomic if it does not contain any atom \cite{Johnson}. Eventually we get the following characterization of non-atomic measure spaces; $(X,\mathcal{A},\mu)$ is non-atomic if and only if every $\mu$-bounded subset of $X$ if of measure zero [Theorem \ref{Th7.10}].
	\section{Some basic facts concerning $U_\mu^I$-topology and $m_\mu^I$-topology}\label{Sec2}
	We begin this section by showing that different ideals $I,J$ in $\mathcal{M}(X,\mathcal{A})$ give rise to different $U_\mu^I, U_\mu^J$-topologies with its analogous counterpart valid for $m_\mu^I, m_\mu^J$-topologies.
	\begin{theorem}\label{Th2.1}
		Let $I\neq J$ be two ideals in $\mathcal{M}(X,\mathcal{A})$. Then $U_\mu^I$-topology $\neq U_\mu^J$-topology (respectively $m_\mu^I$-topology $\neq m_\mu^J$-topology).
	\end{theorem}
	\begin{proof}
		Let $I\neq J$. Choose $g\in I$ such that $g\notin J$. Without loss of generality we can assume that $|g|<1$ on $X$. We assert that $U_\mu(g,J,1)$, which is an open set in the $U_\mu^J$-topology is not open in the $U_\mu^I$-topology. We argue by contradiction and if possible let $U_\mu(g,J,1)$ be open in $U_\mu^I$-topology. Then there exists $\delta>0$ such that $U_\mu(g,I,\delta)\subseteqq U_\mu(g,J,1)$. We observe the that the function $(1+\frac{\delta}{2}).g$ is a member of $U_\mu(g,J,1)$, consequently $\frac{\delta}{2}.g\in J$ -- a contradiction to the initial choice that $g\notin J$ (Mind that $J$ is an ideal in $\mathcal{M}(X,\mathcal{A})$). Then $U_\mu^I$-topology $\neq U_\mu^J$-topology. Finally on using $|g|<1$ on $X$, we observe that for any positive unit $u$ in the ring $\mathcal{M}(X,\mathcal{A})$, $(1+\frac{u}{2}).g$ is a member of $m_\mu(g,I,u)$ without being a member of $m_\mu(g,J,1)$. Therefore $m_\mu(g,J,1)$ is an open set in the $m_\mu^J$-topology without being an open set in the $m_\mu^I$-topology.
	\end{proof}
	However for any chain of ideals in $\mathcal{M}(X,\mathcal{A})$, correspondingly we get chains of $U^I_\mu$ and $m^I_\mu$ topologies on $\mathcal{M}(X,\mathcal{A})$.
	\begin{theorem}
		If $I$ and $J$ are two ideals in $\mathcal{M}(X,\mathcal{A})$ such that $I\subseteqq J$, then the $U^I_\mu$-topology is finer than the $U^J_\mu$-topology on $\mathcal{M}(X,\mathcal{A})$ (and also $m^J_\mu$-topology $\subseteqq m^I_\mu$-topology on $\mathcal{M}(X,\mathcal{A})$).
	\end{theorem}
	The next result characterizes the ideal $\{0\}$ amongst all the ideals in $\mathcal{M}(X,\mathcal{A})$.
	\begin{theorem}
		The following statements are equivalent for an ideal $I$ in $\mathcal{M}(X,\mathcal{A})$: 
		\begin{enumerate}
			\item $I=\{0\}$.
			\item $\mathcal{M}(X,\mathcal{A})$ in the $U_\mu^I$-topology is a discrete space.
			\item $\mathcal{M}(X,\mathcal{A})$ in the $m_\mu^I$-topology is a discrete space.
		\end{enumerate}
	\end{theorem}
	\begin{proof}
		Let $I=\{0\}$. Then for each $g\in \mathcal{M}(X,\mathcal{A})$, $U_\mu(g,I,\frac{1}{2})=\{g\}$. This implies that the $U^I_\mu$-topology is a discrete topology. Since $m^I_\mu$-topology is finer than the $U_\mu^I$-topology, it is also the same as the discrete topology. Thus $(1)\implies(2)$ and $(1)\implies(3)$ hold. To prove the remaining part of the Theorem assume that $I\neq\{0\}$. It follows from Theorem~\ref{Th2.1} that $m_\mu^I$-topology $\neq m_\mu^{\{0\}}$-topology $=$ discrete topology and similarly $U_\mu ^I$-topology $\neq$ discrete topology.
	\end{proof}
	The result that follows decides the nature of the ideals $I$ in $\mathcal{M}(X,\mathcal{A})$ which renders this ring in the $U_\mu^I$-topology a topological ring. 
	\begin{theorem}\label{Th2.3}
		Let $\mathcal{M}(X,\mathcal{A})$ in the $U^I_\mu$-topology be a topological ring, then $I\subseteqq L^\infty(X,\mathcal{A},\mu)$.
	\end{theorem}
	The following Lemma will be helpful towards a proof of Theorem~\ref{Th2.3}.
	\begin{lemma}\label{Lem2.4}
		Let $f\in \mathcal{M}(X,\mathcal{A})\setminus L^\infty(X,\mathcal{A},\mu), f>0$ on $X$. Then there can be constructed a strictly increasing sequence of natural numbers: $1=k_0<k_1<k_2<...$ with the following property: $k_n>n$ for each $n\in\mathbb{N}$ and $\mu\{x\in X: k_i\leq f(x)<k_{i+1}\}>0$ for each $i\in\mathbb{N}\cup\{0\}$. [This Lemma can be proved by using the countable subadditivity of the measure $\mu:\mathcal{A}\to[0,\infty]$, in compliance with the principle of Mathematical induction, see also the proof of $(d)\implies(a)$ in Theorem $3.6$ in~\cite{Acharyya3}].
	\end{lemma}
	\begin{proof}
		of Theorem~\ref{Th2.3}.\\
		If possible let $I\not\subset L^\infty(X,\mathcal{A},\mu)$. Then there exists $f\in I, f>0$ on $X$ such that $f\notin L^\infty(X,\mathcal{A},\mu)$. It follows from Lemma~\ref{Lem2.4} that there exists an increasing sequence of natural numbers $1=k_0<k_1<k_2<...$ with $k_n>n$ for all $n\in\mathbb{N}$ and $\mu(E_i)>0$ for each $i\in\mathbb{N}\cup\{0\}$ where $E_i=\{x\in X:k_i\leq f(x)<k_{i+1}\}$. Since the map:
		\begin{alignat*}{2} \mathcal{M}_{U^I_\mu}(X,\mathcal{A})\times\mathcal{M}_{U^I_\mu}(X,\mathcal{A})&\to\mathcal{M}_{U^I_\mu}(X,\mathcal{A})\\
			(g,h)&\mapsto g.h
		\end{alignat*} is continuous at the point $(0,f)$ by hypothesis [here $\mathcal{M}_{U^I_\mu}(X,\mathcal{A})$ stands for $\mathcal{M}(X,\mathcal{A})$ with the $U^I_\mu$-topology], there exists $\epsilon>0$ such that $U_\mu(0,I,\epsilon)\times U_\mu(f,I,\epsilon)\subseteqq U_\mu(0,I,1)$. Let $g=\frac{\epsilon.f}{2(1+f)}$, then $g\in U_\mu(0,I,\epsilon)$. This implies that $g.f\in U_\mu(0,I,1)$ and hence $\frac{\epsilon.f^2(x)}{2(1+f(x))}<1$ for all $x$ almost everywhere on $X$. This means that there exists $E\in\mathcal{A}$ with $\mu(E)=0$ for which we can write: $\forall x\in X\setminus E,~\frac{\epsilon.f^2(x)}{2(1+f(x))}<1\ldots (A)$. Since $\mu(E_i)>0$ and $\mu(E)=0$, it is clear that $E_i\not\subset E$ for each $i\in\mathbb{N}\cup\{0\}$. So we can pick up a point $x_i\in E_i$ such that $x_i\notin E,~i\in\mathbb{N}\cup\{0\}$. It follows from $(A)$ that for each $n\in\mathbb{N}\cup\{0\},~\frac{\epsilon.f^2(x_n)}{2(1+f(x_n))}<1\ldots (B)$. But $x_n\in E_n\implies f(x_n)\geq k_n>n\implies \lim\limits_{n\to\infty}f(x_n)=\infty\implies \lim\limits_{n\to\infty}\frac{f(x_n)}{1+f(x_n)}=1$ and therefore there exists $k\in\mathbb{N}\cup\{0\}$ such that $\frac{f(x_n)}{1+f(x_n)}>\frac{2}{3}$ for each $n\geq k\implies \forall n\geq k,~\frac{\epsilon.f(x_n)}{2}\frac{2}{3}<\frac{\epsilon.f(x_n)}{2}\frac{f(x_n)}{1+f(x_n)}<1$, from $(B)$. Consequently then $\{f(x_n)\}_n$ becomes a bounded sequence in $\mathbb{R}$ -- a contradiction to the relation $\lim\limits_{n\to\infty}f(x_n)=\infty$, already obtained above.
	\end{proof}
	\begin{theorem}\label{Th2.5}
		For an ideal $I$ in $\mathcal{M}(X,\mathcal{A})$, the following statements are equivalent:
		\begin{enumerate}
			\item $\mathcal{M}_{U^I_\mu}(X,\mathcal{A})$ is a topological vector space. 
			\item $\mathcal{M}_{m^I_\mu}(X,\mathcal{A})(\equiv \mathcal{M}(X,\mathcal{A})$ equipped with $m_\mu^I$-topology$)$ is a topological vector space.
			\item $I=\mathcal{M}(X,\mathcal{A})$ and $\mathcal{M}(X,\mathcal{A})=L^\infty(X,\mathcal{A},\mu)$.
		\end{enumerate}
	\end{theorem}
	\begin{proof}
		First assume that the condition $(3)$ is true. Then $U_\mu^I$-topology $=U_\mu$-topology and $m_\mu^I$-topology $=m_\mu$-topology. It follows from Theorem $3.7$ in~\cite{Acharyya3} [where it is proved amongst other that $L^\infty(X,\mathcal{A},\mu)=\mathcal{M}(X,\mathcal{A})$ if and only if $\mathcal{M}(X,\mathcal{A})$ in the $U_\mu$-topology is a topological vector space] that $(1)$ is true. Furthermore the hypothesis $\mathcal{M}(X,\mathcal{A})=L^\infty(X,\mathcal{A},\mu)$ clearly tells(along with $I=\mathcal{M}(X,\mathcal{A})$) that $U_\mu$-topology $=m_\mu$-topology. Therefore the condition $(2)$ also becomes true. Thus $(3)\implies(1)$ and $(3)\implies(2)$ are established.\\
		Proof of $(1)\implies(3)$: Let $(1)$ be true. Let $f\in\mathcal{M}(X,\mathcal{A})$. Then there exists $\epsilon>0$ in $\mathbb{R}$ such that $(-\epsilon,\epsilon)\times U_\mu(f,I,\epsilon)\subseteqq U_\mu(0,I,1)$ -- this follows on using the continuity of the scalar multiplication function at the point $(0,f)$. This implies that $\frac{\epsilon}{2}f\in I$ and hence $f\in I$. Thus $\mathcal{M}(X,\mathcal{A})=I$ and clearly then $U_\mu^I$-topology $=U_\mu$-topology. Since now $U_\mu$-topology makes $\mathcal{M}(X,\mathcal{A})$ a topological vector space on account of the hypothesis $(1)$, it follows on using Theorem $3.7$ in~\cite{Acharyya3} once again that $L^\infty(X,\mathcal{A},\mu)=\mathcal{M}(X,\mathcal{A})$. Thus $(1)\implies(3)$ is proved and hence $(1)\iff(3)$.\\
		$(2)\implies(3)$: Let $(2)$ be true. Then by closely following the chain of arguments in the proof of first part of the implication relation $(1)\implies(3)$ above and making some necessary modification, it is not hard to show that $\mathcal{M}(X,\mathcal{A})=I$. Therefore $m_\mu^I$-topology $=m_\mu$-topology and we can say by $(2)$ that $\mathcal{M}(X,\mathcal{A})$ becomes a topological vector space in the $m_\mu$-topology. We assert that $\mathcal{M}(X,\mathcal{A})=L^\infty(X,\mathcal{A},\mu)$ and that finishes the proof.\\
		Proof of the assertion: If possible let $L^\infty(X,\mathcal{A},\mu)\subsetneqq\mathcal{M}(X,\mathcal{A})$. Then there exists $f\in \mathcal{M}(X,\mathcal{A}),f>0$ on $X$ such that $f\notin L^\infty(X,\mathcal{A},\mu)$. Then for each $n\in\mathbb{N}$, $\mu\{x\in X:f(x)>n\}>0$. Take $h=\frac{1}{f}$, then this implies for each $n\in\mathbb{N}$, $\mu\{x\in X:h(x)<\frac{1}{n}\}>0\ldots(A)$. We now show that the scalar multiplication map \begin{alignat*}{1}
			S:\mathbb{R}\times\mathcal{M}(X,\mathcal{A})&\to\mathcal{M}(X,\mathcal{A})\\
			(\lambda,g)&\mapsto\lambda.g
		\end{alignat*} is not continuous at the point $(1,\underline{1})$, here $\underline{1}$ is the constant function on $X$ with value $1$ and that gives the desired contradiction. If possible let $S$ be continuous at $(1,\underline{1})$. Then there exists $\alpha>0$ in $\mathbb{R}$ and a positive unit $u$ in the ring $\mathcal{M}(X,\mathcal{A})$ for which we can write: $(1-\alpha,1+\alpha)\times m_\mu(\underline{1},u)\subseteqq m_\mu(\underline{1},h)$ [We write for simplicity for this choice of $I=\mathcal{M}(X,\mathcal{A})$, $m_\mu(\underline{1},u)=m_\mu(\underline{1},I,u)$ etc.]. This implies that $|\underline{(1+\frac{\alpha}{2})}-\underline{1}|<h$ almost everywhere on $X$, i.e., $h>\frac{\alpha}{2}$ almost everywhere on $X$. Now there exists $n\in\mathbb{N}$ such that $\frac{\alpha}{2}>\frac{1}{n}$ and hence $h>\frac{1}{n}$ almost everywhere on $X$. Consequently $\mu\{x\in X:h(x)\leq\frac{1}{n}\}=0$ -- this contradicts the relation $(A)$ obtained above. The proof is now complete.
	\end{proof}
	The following condition on the nature of the ideals $I$, amongst a class of ideals in $\mathcal{M}(X,\mathcal{A})$, for the coincidence of $U_\mu^I$-topology and $m_\mu^I$-topology on $\mathcal{M}(X,\mathcal{A})$ seems to be a fascinating one.
	\begin{theorem}\label{Th2.6}
		For any ideal $I$ in $\mathcal{M}(X,\mathcal{A})$ if $X\setminus\bigcap Z[I]$ is a $\mu$-bounded subset of $X$, then $U_\mu^I$-topology $=m_\mu^I$-topology on $\mathcal{M}(X,\mathcal{A})$. Moreover if $I$ is countably generated then the converse is also true.
	\end{theorem}
	\begin{proof}
		First assume that $X\setminus\bigcap Z[I]$ is a $\mu$-bounded subset of $X$. A typical basic open set in the $m_\mu^I$-topology is of the form $m_\mu(f,I,u)$, here $f\in\mathcal{M}(X,\mathcal{A})$ and $u$, a positive unit in the latter ring. The hypothesis tells that there exists $\lambda>0$ such that $\frac{1}{u(x)}<\lambda$ holds for all $x$, almost everywhere on $X\setminus\bigcap Z[I]$. This means that there exists a subset $E_0$ of $X\setminus\bigcap Z[I]$ with $\mu(E_0)=0$ for which for each $x\in (X\setminus\bigcap Z[I])\setminus E_0,~u(x)>\frac{1}{\lambda}$. Now let $g\in U_\mu(f,I,\frac{1}{\lambda})$, then there exists a measurable set $A_g$ with $\mu(A_g)=0$ for which we can write: $\sup\limits_{x\in X\setminus A_g}|g(x)-f(x)|<\frac{1}{\lambda}$ and also we have $g-f\in I$. This implies that for each $x\in X\setminus(A_g\cup E_0),~|f(x)-g(x)|<u(x)$ and consequently then $g\in m_\mu(f,I,u)$. Thus it is proved that $U_\mu(f,I,\frac{1}{\lambda})\subseteqq m_\mu(f,I,u)$. Hence $m_\mu(f,I,u)$ is open in the $U_\mu^I$-topology and consequently, $m_\mu^I$-topology $\subseteqq U_\mu^I$-topology $\subseteqq m_\mu^I$-topology. Therefore $U_\mu^I$-topology $=m_\mu^I$-topology.\\
		To prove the converse assume that $U_\mu^I$-topology $=m_\mu^I$-topology; now assume that $I$ is countably generated. If possible let $Y=X\setminus\bigcap Z[I]$ be not $\mu$-bounded. Then there exists $f\in\mathcal{M}(X,\mathcal{A})$, such that $f>1$ on $X$ and $f$ is not essentially bounded on $Y$. On using the result of Lemma~\ref{Lem2.4}, we can produce a strictly increasing sequence of natural numbers: $1=k_0<k_1<k_2<...$ with the following property: $k_n>n$ for each $n\in\mathbb{N}$ and for each $i\in\mathbb{N}\cup\{0\},~\mu(E_i)>0$, where $E_i=\{y\in Y:k_i\leq f(y)<k_{i+1}\}$. Then $Y=\bigcup\limits_{n=0}^\infty E_n=$ a countable union of pairwise disjoint measurable sets with each $\mu(E_n)>0$. Define $u:X\to\mathbb{R}$ as follows: $u(x)=\begin{cases}
			\frac{1}{f(x)}&\text{ if }x\in Y\\
			1&\text{ if }x\in X\setminus Y
		\end{cases}$. Then $u\in\mathcal{M}(X,\mathcal{A})$ and is a positive unit in this ring and we can now write: $E_i=\{y\in Y:\frac{1}{k_{i+1}}<u(y)\leq\frac{1}{k_i}\}\ldots(A)$. We shall show that $m_\mu(0,I,u)$, which is an open set in the $m_\mu^I$-topology is not open in the $U_\mu^I$-topology and that gives the necessary contradiction.\\
	 	Towards proving this result we first claim the followings: given $i\in\mathbb{N}$, there exists $j\geq i$ and $g\in I$ such that $g\neq 0$ throughout a set of positive $\mu$-measure contained in $E_g\ldots(B)$.\\
		If possible let the claim $(B)$ made above be false. Choose $g\in I$ and $j\geq i$, arbitrarily. Now since $I$ is countably generated, say $I=<g_n>_{n=1}^\infty$, it follows that for each $n\in\mathbb{N}, ~g_n=0$ almost everywhere on $E_j$. This means that there exists a measurable set $F_{n,j}\subseteqq E_j$ such that $\mu(F_{n,j})=0$ and $g_n=0$ on $E_j\setminus F_{n,j}$. Let $F_j=\bigcup\limits_{n=1}^\infty F_{n,j}$, then $F_j\subseteqq E_j$ and $\mu(F_j)=0$. Also we note that for the choice $g\in I,~g=0$ on $E_j\setminus F_j$, because $g$ is a finite linear combination of $g_n$'s with coefficients from $\mathcal{M}(X,\mathcal{A})$. Thus we can say that for each $j\geq i$, there exists a measurable set $F_J\subseteqq E_j$ with $\mu(F_j)=0$ for which $\forall g\in I,~g=0$ on $E_j\setminus F_j$. This further implies that $\bigcup\limits_{j=1}^\infty (E_j\setminus F_j)=Y_0$ (say) $\subseteqq \bigcap Z[I]$. On the other hand $Y_0\subseteqq Y=X\setminus\bigcap Z[I]$. These last two statements are consistent provided that $Y_0=\emptyset$, which is never the case because for each $j\geq i,~\mu(E_j\setminus F_j)>0\implies E_j\setminus F_j\neq\emptyset~\forall j\geq i\implies Y_0\neq\emptyset$. This contradiction proves that the claim $(B)$ made above is right. We recall that in order to complete the proof we should ascertain that $m_\mu(0,I,u)$ is not open in the $U_\mu^I$-topology. If possible let $m_\mu(0,I,u)$ be open in the $U_\mu^I$-topology. Then there exists $\epsilon>0$ in $\mathbb{R}$ such that $U_\mu(0,I,\epsilon)\subseteqq m_\mu(0,I,u)\ldots(C)$. Since $\lim\limits_{n\to\infty}k_n=\infty$, there exists $i\in\mathbb{N}$ such that for all $j\geq i,~\frac{2}{k_j}<\frac{\epsilon}{2}$. This implies in view of $(C)$ that $U_\mu(0,I,\frac{2}{k_j})\subseteqq m_\mu(0,I,u)$ for each $j\geq i\ldots(D)$. Now in view of the claim $(B)$ made above; there exists $j\geq i,~g\in I$ and $F_j\subseteqq E_j$ with $\mu(E_j\setminus F_j)>0$ such that $g\neq 0$ throughout the whole of the set $E_j\setminus F_j$. Define $h:X\to\mathbb{R}$ as follows: $h(x)=\begin{cases}
			\frac{1}{k_jg(x)}&\text{ if }x\in F_j\\
			1&\text{ if }x\in X\setminus F_j
		\end{cases}$. Then $h\in\mathcal{M}(X,\mathcal{A})$ and hence $g.h\in I$. Let $l=gh\wedge\frac{1}{k_j}\in\mathcal{M}(X,\mathcal{A})$. Since each ideal in $\mathcal{M}(X,\mathcal{A})$ is absolutely convex [Proposition $3.2$,~\cite{Acharyya1}], it follows that $l\in I$ and of course $l\in U_\mu(0,I,\frac{2}{k_j})$. It follows from $(D)$therefore that $l\in m_\mu(0,I,u)$. Consequently $l(x)<u(x)$ almost everywhere on $X$ and therefore $\mu\{x\in X: l(x)\geq u(x)\}=0\ldots(E)$. But since $l(x)=\frac{1}{k_j}$ on $F_j$ and $u(x)\leq\frac{1}{k_j}$ on $F_j$ (because $F_j\subseteqq E_j$ and $u\leq\frac{1}{k_j}$ on $E_j$ from $(A)$), it follows that $l(x)\geq u(x)$ on $F_j$ with $\mu(F_j)>0$. This implies that $\mu\{x\in X: l(x)\geq u(x)\}>0$ -- a contradiction to the relation $(E)$.
	\end{proof}
	\begin{corollary}\label{Cor2.7}
		$U_\mu$-topology $=m_\mu$-topology if and only if each function in $\mathcal{M}(X,\mathcal{A})$ is essentially bounded, i.e., $\mathcal{M}(X,\mathcal{A})=L^\infty(X,\mathcal{A},\mu)$.
	\end{corollary}
	[The proof follows on choosing $I=\mathcal{M}(X,\mathcal{A})$.]
	\section{The component of $0$ in the $U_\mu^I$-topology and $m_\mu^I$-topology}\label{Sec3}
	At the very outset in this section we make it clear that we let $\mathcal{M}_{U_\mu^I}(X,\mathcal{A})$ to stand for $\mathcal{M}(X,\mathcal{A})$ equipped with $U_\mu^I$-topology with an analogous meaning for $\mathcal{M}_{m_\mu^I}(X,\mathcal{A})$. We start with finding out two special clopen subsets of $\mathcal{M}_{U_\mu^I}(X,\mathcal{A})$, the first one containing $I$ and the second contained in $I$.
	\begin{theorem}\label{Th3.1}
		\hspace*{3cm}
		\begin{enumerate}
			\item If $J$ is any additive subgroup of $\mathcal{M}(X,\mathcal{A})$ containing $I$, then $J$ is clopen in $\mathcal{M}_{U_\mu^I}(X,\mathcal{A})$.
			\item $I\cap L^\infty(X,\mathcal{A},\mu)$ is a clopen subset of $\mathcal{M}_{U_\mu^I}(X,\mathcal{A})$.\label{Th3.1(2)}
		\end{enumerate}
	\end{theorem}
	\begin{proof}
		\hspace*{3cm}
		\begin{enumerate}
			\item For any $f\in J$, it is easy to prove that $f\in U_\mu(f,I,1)\subset J$. This proves that $J$ is open in $\mathcal{M}_{U_\mu^I}(X,\mathcal{A})$. On the other hand for any $f\in\mathcal{M}(X,\mathcal{A})\setminus J,~U_\mu(f,I,1)\cap J=\emptyset$ -- an easy check. This proves that $J$ is closed in $\mathcal{M}_{U_\mu^I}(X,\mathcal{A})$.
			\item It is easy to check that for any $f\in I\cap L^\infty(X,\mathcal{A},\mu),~f\in U_\mu(f,I,1)\subset I\cap L^\infty(X,\mathcal{A},\mu)$. Hence $I\cap L^\infty(X,\mathcal{A},\mu)$ is open in $\mathcal{M}_{U_\mu^I}(X,\mathcal{A})$. On the other hand if $f\in\mathcal{M}(X,\mathcal{A})\setminus(I\cap L^\infty(X,\mathcal{A},\mu))$, then it is not hard to check that $U_\mu(f,I,1)\cap I\cap L^\infty(X,\mathcal{A},\mu)=\emptyset$. Thus $I\cap L^\infty(X,\mathcal{A},\mu)$ is clopen in $\mathcal{M}_{U_\mu^I}(X,\mathcal{A})$.
		\end{enumerate}
	\end{proof}
	By borrowing idea from~\cite{Azarpanah5} [Lemma $3.1$], for any $f\in \mathcal{M}(X,\mathcal{A})$, we define the map \begin{alignat*}{1}
		\phi_f:\mathbb{R}&\to\mathcal{M}(X,\mathcal{A})\\
		r&\mapsto r.f
	\end{alignat*}
	\begin{theorem}\label{Th3.2}
		\hspace*{3cm}
		\begin{enumerate}
			\item For each $f\in\mathcal{M}(X,\mathcal{A}),~\phi_f:\mathbb{R}\to\mathcal{M}_{U_\mu^I}(X,\mathcal{A})$ is a continuous map if and only if $f\in I\cap L^\infty(X,\mathcal{A},\mu)$.\label{Th3.2(1)}
			\item $I\cap L^\infty(X,\mathcal{A},\mu)$ is the component of $0$ in $\mathcal{M}_{U_\mu^I}(X,\mathcal{A})$.
		\end{enumerate}
	\end{theorem}
	\begin{proof}
		\hspace*{3cm}
		\begin{enumerate}
			\item Let $\phi_f$ be a continuous map, in particular at the point $r=0$. Then there exists $\delta>0$ in $\mathbb{R}$ such that $\phi_f(-\delta,\delta)\subseteqq U_\mu(\phi_f(0),I,1)=U_\mu(0,I,1)$. It follows that $\phi_f(\frac{\delta}{2})\in I$. This implies that $\frac{\delta}{2}.f\in I$ and also $|\frac{\delta}{2}.f|<1$ almost everywhere on $X$. Hence $f\in I$ and $f$ is essentially bounded on $X$, thus $f\in I\cap L^\infty(X,\mathcal{A},\mu)$.\\
			Conversely let $f\in I\cap L^\infty(X,\mathcal{A},\mu)$. Then $|f|\leq\lambda$ almost everywhere on $X$ for some $\lambda>0$. Choose $r\in\mathbb{R}$ and $\epsilon>0$ at random. Then $\phi_f(r-\frac{\epsilon}{\lambda}, r+\frac{\epsilon}{\lambda})\subset U_\mu(\phi_f(r),I,\epsilon)$. This ensures the continuity of $\phi_f$ at $r$.
			\item We see that $\bigcup\limits_{f\in I\cap L^\infty(X,\mathcal{A},\mu)}\phi_f(\mathbb{R})=I\cap L^\infty(X,\mathcal{A},\mu)\equiv$ the union of a family of connected subsets of $\mathcal{M}_{U_\mu^I}(X,\mathcal{A})$, each containing the function $0$. Hence $I\cap L^\infty(X,\mathcal{A},\mu)$ is a connected subset of the latter space. Since by Theorem~\ref{Th3.1}(\ref{Th3.1(2)}), this last set is already clopen in $\mathcal{M}_{U_\mu^I}(X,\mathcal{A})$, it follows that, $I\cap L^\infty(X,\mathcal{A},\mu)$ is the component of $0$ in $\mathcal{M}_{U_\mu^I}(X,\mathcal{A})$.
		\end{enumerate}
	\end{proof}
	\begin{corollary}
		The component of $0$ in $\mathcal{M}(X,\mathcal{A})$ in the $U_\mu$-topology is $L^\infty(X,\mathcal{A},\mu)$. [This follows on choosing $I=\mathcal{M}(X,\mathcal{A})$]
	\end{corollary}
	The following subset will be turn out to be crucial for finding out the component of $0$ in $\mathcal{M}_{m_\mu^I}(X,\mathcal{A})\equiv$ the set $\mathcal{M}(X,\mathcal{A})$ equipped with $m_\mu^I$-topology.
	\begin{notation}
		We set $L_\psi(X,\mathcal{A},\mu)=\{f\in\mathcal{M}(X,\mathcal{A}):\forall g\in\mathcal{M}(X,\mathcal{A}),f.g\in L^\infty(X,\mathcal{A},\mu)\}$.
	\end{notation}
	It is easy to verify that $L_\psi(X,\mathcal{A},\mu)$ is an ideal in the ring $\mathcal{M}(X,\mathcal{A})$ and is $\mu$-bounded in the sense that $L_\psi(X,\mathcal{A},\mu)\subseteqq L^\infty(X,\mathcal{A},\mu)$. Furthermore it is not at all difficult to check that $L_\psi(X,\mathcal{A},\mu)$ is the largest $\mu$-bounded ideal in $\mathcal{M}(X,\mathcal{A})$: indeed if $J$ is an ideal in $\mathcal{M}(X,\mathcal{A})$ and $J\subseteqq L^\infty(X,\mathcal{A},\mu)$, then $J\subseteqq L_\psi(X,\mathcal{A},\mu)$ for if there exists $f\in J$ such that $f\notin L_\psi(X,\mathcal{A},\mu)$, then $f.g\notin L^\infty(X,\mathcal{A},\mu)$ for some $g\in\mathcal{M}(X,\mathcal{A})$; but $f\in J$ and $J$ is an ideal in $\mathcal{M}(X,\mathcal{A})\implies f.g\in J\implies f.g\in L^\infty(X,\mathcal{A},\mu)$, a contradiction. We now offer an alternative description of $L_\psi(X,\mathcal{A},\mu)$.
	\begin{theorem}\label{Th3.55}
		$L_\psi(X,\mathcal{A},\mu)=\{f\in\mathcal{M}(X,\mathcal{A}):X\setminus Z(f)\text{ is }\mu\text{-bounded}\}$
	\end{theorem}
	\begin{proof}
		Let $X\setminus Z(f)$ be $\mu$-bounded for some $f\in \mathcal{M}(X,\mathcal{A})$ and $g\in\mathcal{M}(X,\mathcal{A})$, arbitrary. Then there exists  $\mu$-null sets $E_1, E_2\subseteqq X\setminus Z(f)$ such that $|g|<\lambda$ on $(X\setminus Z(f))\setminus E_1$ and $|f|<\delta$ on $(X\setminus Z(f))\setminus E_2$ for some $\lambda,\delta>0$. Therefore $|f.g|<\lambda.\delta$ on $X\setminus (E_1\cup E_2)\implies f\in L_\psi(X,\mathcal{A},\mu)$. Conversely let $f\in L_\psi(X,\mathcal{A},\mu)$. Define $h:X\to\mathbb{R}$ as follows: $h(x)=\begin{cases}
			\frac{1}{f(x)}&\text{ if }x\notin Z(f)\\
			1&\text{ if }x\in Z(f)
		\end{cases}$. Then $h\in \mathcal{M}(X,\mathcal{A})$. Let $g\in\mathcal{M}(X,\mathcal{A})$, arbitrary. Then $f\in L_\psi(X,\mathcal{A},\mu)\implies f.g.h\in L^\infty(X,\mathcal{A},\mu)$. So there exists $E\in\mathcal{A}$ with $\mu(E)=0$ such that $|f.g.h|<\lambda$ on $X\setminus E$ for some $\lambda>0$. Thus for all $x\in (X\setminus Z(f))\setminus E$, $f(x).h(x)=1\implies |g(x)|<\lambda$. Let $E_1=X\setminus Z(f)\cap E$. Then $E_1\in\mathcal{A}$ and $E_1$ is a $\mu$-null subset of $X\setminus Z(f)$ and $g$ is bounded on $(X\setminus Z(f))\setminus E_1$. Since $g\in\mathcal{M}(X,\mathcal{A})$ is arbitrary, it follows that $X\setminus Z(f)$ is $\mu$-bounded.
	\end{proof}
	\begin{remark}
		$L_\psi(X,\mathcal{A},\mu)$ is the measure theoretic analogue of the well-known ideal $C_\psi(X)$ of $C(X)$, here $C_\psi(X)$ is the set of all those functions in $C(X)$ which have pseudocompact support. [See Theorem $2.1$,\cite{Mandelker}]
	\end{remark}
	The next result decides the component of $0$ in $\mathcal{M}_{m_\mu^I}(X,\mathcal{A})$ in clear terms.
	\begin{theorem}\label{Th3.5}
		The component of $0$ in $\mathcal{M}_{m_\mu^I}(X,\mathcal{A})$ is $I\cap L_\psi(X,\mathcal{A},\mu)$.
	\end{theorem}
	We need the following Lemma to prove this fact:
	\begin{lemma}
		$f\in I\cap L_\psi(X,\mathcal{A},\mu)$ if and only if the given map is continuous \begin{alignat*}{1}
			\phi_f:\mathbb{R}&\to\mathcal{M}_{m_\mu^I}(X,\mathcal{A})\\
			r&\mapsto r.f
		\end{alignat*}
	\end{lemma}
	\begin{proof}
		First assume that $f\in I\cap L_\psi(X,\mathcal{A},\mu)$. Let $u$ be a positive unit in $\mathcal{M}(X,\mathcal{A})$. Define $g:X\to\mathbb{R}$ as follows: $g(x)=\begin{cases}
			\frac{1}{f(x)}&\text{ if }x\in X\setminus Z(f)\\
			1&\text{ if }x\in Z(f)
		\end{cases}$. Then $g\in\mathcal{M}(X,\mathcal{A})$. As $f\in L_\psi(X,\mathcal{A},\mu)$, this implies that $\frac{g.f}{u}\in L^\infty(X,\mathcal{A},\mu)$. This further implies that $\frac{1}{u}$ is essentially bounded on $X\setminus Z(f)$, as $g.f=1$ on $X\setminus Z(f)$. Thus $u\geq\frac{1}{\lambda}$ almost everywhere on $X\setminus Z(f)$ for some $\lambda>0$. Now for $r\in\mathbb{R},~m_\mu(r.f,I,u)$ is a basic open neighbourhood of $\phi_f(r)=r.f$ in $\mathcal{M}_{m_\mu^I}(X,\mathcal{A})$. As $f\in L_\psi(X,\mathcal{A},\mu)\subseteqq L^\infty(X,\mathcal{A},\mu),~|f|\leq M$ for some $M>0$ almost everywhere on $X$. Let $s\in\mathbb{R}$ be such that $|s-r|<\frac{1}{\lambda.M}$. Then $|\phi_f(s)-\phi_f(r)|=|s.f-r.f|<\frac{1}{\lambda.M}.M$ almost everywhere on $X$. This implies that $|s.f-r.f|<u$ almost everywhere on $X$, i.e., $s.f\in m_\mu(r.f,I,u)$ [Note that $f\in I\implies s.f-r.f\in I$]. Thus $\phi_f$ becomes continuous at $r\in\mathbb{R}$.\\
		Conversely let \begin{alignat*}{1}
			\phi_f:\mathbb{R}&\to\mathcal{M}_{m_\mu^I}(X,\mathcal{A})\\
			r&\mapsto r.f
		\end{alignat*} be continuous for some $f\in\mathcal{M}(X,\mathcal{A})$. We shall show that $f\in I\cap L_\psi(X,\mathcal{A},\mu)$. Now $U_\mu^I$-topology $\subseteqq m_\mu^I$-topology $\implies$ the map \begin{alignat*}{1}
		\phi_f:\mathbb{R}&\to\mathcal{M}_{U_\mu^I}(X,\mathcal{A})\\
		r&\mapsto r.f
		\end{alignat*} is continuous. It follows from 	Theorem~\ref{Th3.2}(\ref{Th3.2(1)}) that $f\in I\cap L^\infty(X,\mathcal{A},\mu)$. To show that $f\in L_\psi(X,\mathcal{A},\mu)$, choose $g\in\mathcal{M}(X,\mathcal{A})$. We need to verify that $f.g\in L^\infty(X,\mathcal{A},\mu)$. Take $u=\frac{1}{1+|g|}$. Since \begin{alignat*}{2}
			\phi_f:\mathbb{R}&\to\mathcal{M}_{m_\mu^I}(X,\mathcal{A})\\
			r&\mapsto r.f
		\end{alignat*} is continuous at $r=0$, there exists $\delta>0$ in 	$\mathbb{R}$ such that $\phi_f(-\delta,\delta)\subseteqq m_\mu(\phi_f(0),I,u)=m_\mu(0,I,u)$. So $\phi_f(\frac{\delta}{2})\in m_\mu(0,I,u)\implies |\frac{\delta}{2}f|<u=\frac{1}{1+|g|}$ almost everywhere on $X$. Consequently $f.g\in L^\infty(X,\mathcal{A},\mu)$. The lemma is proved.
	\end{proof}
	\begin{proof} of the main Theorem:
		It follows form the Lemma that $I\cap L_\psi(X,\mathcal{A},\mu)=\bigcup\limits_{f\in I\cap L_\psi(X,\mathcal{A},\mu)}\phi_f(\mathbb{R})\equiv$ a connected set in $\mathcal{M}_{m_\mu^I}(X,\mathcal{A})$. Thus $I\cap L_\psi(X,\mathcal{A},\mu)$ is a connected ideal in $\mathcal{M}_{m_\mu^I}(X,\mathcal{A})$. To complete this Theorem we need to prove that $I\cap L_\psi(X,\mathcal{A},\mu)$ is a maximal connected ideal in $\mathcal{M}_{m_\mu^I}(X,\mathcal{A})$. For that purpose let $J$ be any connected ideal in $\mathcal{M}_{m_\mu^I}(X,\mathcal{A})$. Now by Theorem~\ref{Th3.1}(\ref{Th3.1(2)}), $I\cap L^\infty(X,\mathcal{A},\mu)$ is a clopen subset of $\mathcal{M}_{U_\mu^I}(X,\mathcal{A})$. Since the $m_\mu^I$-topology is finer than the $U_\mu^I$-topology on $\mathcal{M}(X,\mathcal{A})$, it follows that $I\cap L^\infty(X,\mathcal{A},\mu)$ is a clopen subset in the space $\mathcal{M}_{m_\mu^I}(X,\mathcal{A})$. Since $J$ is a connected ideal in $\mathcal{M}_{m_\mu^I}(X,\mathcal{A})$, this implies that $J\subseteqq I\cap L^\infty(X,\mathcal{A},\mu)$. In particular $J$ becomes a $\mu$-bounded ideal in $\mathcal{M}(X,\mathcal{A})$. As $L_\psi(X,\mathcal{A},\mu)$ is the largest $\mu$-bounded ideal in $\mathcal{M}(X,\mathcal{A})$, already observed, it follows that $J\subseteqq L_\psi(X,\mathcal{A},\mu)$. Consequently, $J\subseteqq I\cap L_\psi(X,\mathcal{A},\mu)$.
	\end{proof}
	\begin{corollary}
		The component of $0$ in $\mathcal{M}(X,\mathcal{A})$ equipped with the $m_\mu$-topology is $L_\psi(X,\mathcal{A},\mu)$.
	\end{corollary}
	\section{Closures of chosen ideals of $\mathcal{M}(X,\mathcal{A})$ in the $U_\mu\backslash m_\mu$-topology}\label{Sec4}
	In this section and also in section \ref{Sec5}, we assume throughout that $I=\mathcal{M}(X,\mathcal{A})$. Therefore $U_\mu^I$-topology $=U_\mu$-topology and $m_\mu^I$-topology $=m_\mu$-topology. It is a standard result in the theory of $C(X)$, the ring of real valued continuous functions defined over a Tychonoff space $X$ that $C_\infty(X)=$ the closure of $C_k(X)$ in $C(X)$ in the topology of uniform convergence. Here $C_\infty(X)=\{f\in C(X):\forall\epsilon>0,\{x\in X:|f(x)|\geq\epsilon\}\text{ is compact}\}$ and $C_k(X)=\{f\in C_\infty(X):cl_X(X\setminus Z(f))\text{ is compact}\}$. In what follows we determine the appropriate measure theoretic analogue of this fact.
	\begin{theorem}\label{Th4.0}
		$L^\infty(X,\mathcal{A},\mu)$ is a closed set in $\mathcal{M}(X,\mathcal{A})$ in the $U_\mu$-topology.
	\end{theorem}
	\begin{proof}
		Let $f\in\mathcal{M}(X,\mathcal{A})\setminus L^\infty(X,\mathcal{A},\mu)$. If possible let $U_\mu(f,1)\cap L^\infty(X,\mathcal{A},\mu)\neq\emptyset$. Then there exists $g\in L^\infty(X,\mathcal{A},\mu)$ such that $g\in U_\mu(f,1)$. So there exist $A_1,A_2\in\mathcal{A}$ with $\mu(A_1)=0=\mu(A_2)$ such that $|f-g|<1$ on $X\setminus A_1$ and $|g|<\lambda$ on $X\setminus A_2$ for some $\lambda>0$. Therefore $|f|\leq |f-g|+|g|<1+\lambda$ on $X\setminus(A_1\cup A_2)\implies f\in L^\infty(X,\mathcal{A},\mu)$ -- a contradiction. Therefore $U_\mu(f,1)\cap L^\infty(X,\mathcal{A},\mu)=\emptyset$. Consequently, $L^\infty(X,\mathcal{A},\mu)$ is closed in the $U_\mu$-topology.
	\end{proof}
	\begin{theorem}\label{Th4.1}
		$M_\infty(X,\mathcal{A})=cl_{U_\mu}M_F(X,\mathcal{A})\equiv$ the closure of $M_F(X,\mathcal{A})$ in the $U_\mu$-topology.
	\end{theorem}
	\begin{proof}
		We first show that $M_\infty(X,\mathcal{A})$ is closed in the space $\mathcal{M}_{U_\mu}(X,\mathcal{A})\equiv \mathcal{M}(X,\mathcal{A})$ equipped with the $U_\mu$-topology. So let $f\in\mathcal{M}(X,\mathcal{A})\setminus\mathcal{M}_\infty(X,\mathcal{A})$. Then there exists $n\in\mathbb{N}$ such that $\mu\{x\in X:|f(x)|\geq\frac{1}{n}\}=\infty\ldots(A)$. We claim that $U_\mu(f,\frac{1}{2n})\cap\mathcal{M}_\infty(X,\mathcal{A})=\emptyset$ and that will do. Here $U_\mu(f,\frac{1}{2n})=\{g\in \mathcal{M}(X,\mathcal{A}):\sup\limits_{x\in X\setminus A_g} |f(x)-g(x)|<\frac{1}{2n}\text{ for some }A_g\in\mathcal{A}\text{ with }\mu(A_g)=0\}$. If possible let $U_\mu(f,\frac{1}{2n})\cap\mathcal{M}_\infty(X,\mathcal{A})\neq\emptyset$, choose $g$ from this non-empty set. Then there exists a measurable set $A_g$ with $\mu(A_g)=0$ such that $\sup\limits_{x\in X\setminus A_g}|f(x)-g(x)|<\frac{1}{2n}$. But $g\in\mathcal{M}_\infty(X,\mathcal{A})\implies\mu\{x\in X:|g(x)|\geq\frac{1}{2n}\}<\infty$. Since $\{x\in X:|f(x)|\geq\frac{1}{n}\}\subseteqq\{x\in X:|g(x)|\geq\frac{1}{2n}\cup\{x\in X:|f(x)-g(x)|\geq\frac{1}{2n}\}\subseteqq\{x\in X:|g(x)|\geq\frac{1}{2n}\}\cup A_g$, it follows that $\mu\{x\in X:|f(x)|\geq\frac{1}{n}\}\leq\mu\{x\in X:|g(x)|\geq\frac{1}{2n}\}$ (as $\mu(A_g)=0$) and consequently then $\mu\{x\in X:|f(x)|\geq\frac{1}{n}\}<\infty$, a contradiction to the relation $(A)$. Thus it is realized that $\mathcal{M}_\infty(X,\mathcal{A})$ is closed in $\mathcal{M}_{U_\mu}(X,\mathcal{A})$.\\
		To complete this Theorem let $f\in \mathcal{M}_\infty(X,\mathcal{A})$ and $\epsilon>0$ be preassigned. It is sufficient to show that $U_\mu(f,\epsilon)\cap\mathcal{M}_F(X,\mathcal{A})\neq\emptyset$. Towards that end define a function $g:X\to\mathbb{R}$ as follows: $$g(x)=\begin{cases}
			f(x)+\frac{\epsilon}{2}&\text{ if }f(x)\leq -\frac{\epsilon}{2}\\
			0&\text{ if }-\frac{\epsilon}{2}\leq f(x)\leq\frac{\epsilon}{2}\\
			f(x)-\frac{\epsilon}{2}&\text{ if }f(x)\geq\frac{\epsilon}{2}
		\end{cases}$$Then $g\in\mathcal{M}(X,\mathcal{A})$ and $|f(x)-g(x)|\leq\frac{\epsilon}{2}$ for all $x\in X$. This implies that $g\in U_\mu(f,\epsilon)$. We further assert that $g\in\mathcal{M}_F(X,\mathcal{A})$.\\
		Proof of the assertion: $X\setminus Z(g)\subseteqq\{x\in X:|f(x)|\geq\frac{\epsilon}{2}\}\implies\mu(X\setminus Z(g))\leq\mu\{x\in X:|f(x)|\geq\frac{\epsilon}{2}\}<\infty$ because $f\in\mathcal{M}_\infty(X,\mathcal{A})$. It follows that $g\in \mathcal{M}_F(X,\mathcal{A})$ and hence $g\in U_\mu(f,\epsilon)\cap \mathcal{M}_F(X,\mathcal{A})$.
	\end{proof}
	\begin{remark}
		Since the $m_\mu$-topology on $\mathcal{M}(X,\mathcal{A})$ is finer than the $U_\mu$-topology, it follows from the above Theorem that $cl_{m_\mu}\mathcal{M}_F(X,\mathcal{A})\equiv$ the closure of $\mathcal{M}_F(X,\mathcal{A})$ in the $m_\mu$-topology is contained in $\mathcal{M}_\infty(X,\mathcal{A})$. Indeed the next proposition exactly determines this closure.
	\end{remark}
	\begin{notation}
		Let $M_\psi(X,\mathcal{A})=\{f\in\mathcal{M}(X,\mathcal{A}):\forall g\in\mathcal{M}(X,\mathcal{A}), f.g\in \mathcal{M}_\infty(X,\mathcal{A})\}$.
	\end{notation}
	Then clearly $\mathcal{M}_\psi(X,\mathcal{A})\subseteqq\mathcal{M}_\infty(X,\mathcal{A})$.
	\begin{theorem}\label{Th4.4}
		$\mathcal{M}_\psi(X,\mathcal{A})=cl_{m_\mu}\mathcal{M}_F(X,\mathcal{A})$.
	\end{theorem}
	\begin{proof}
		We first show that $\mathcal{M}_\psi(X,\mathcal{A})$ is closed in $\mathcal{M}_{m_\mu}(X,\mathcal{A})\equiv \mathcal{M}(X,\mathcal{A})$ equipped with the $m_\mu$-topology. Towards that end let for $g\in \mathcal{M}(X,\mathcal{A}),~T_g=\{f\in\mathcal{M}(X,\mathcal{A}):f.g\in\mathcal{M}_\infty(X,\mathcal{A})\}$. Then $\mathcal{M}_\psi(X,\mathcal{A})=\bigcap \{T_g:g\in\mathcal{M}(X,\mathcal{A})\}$. It is therefore sufficient to prove for an arbitrary selected $g\in\mathcal{M}(X,\mathcal{A})$ that $T_g$ is closed in $\mathcal{M}_{m_\mu}(X,\mathcal{A})$. So we choose any function $h\in\mathcal{M}(X,\mathcal{A})\setminus T_g$. Then $h.g\notin\mathcal{M}_\infty(X,\mathcal{A})$ and hence $\mu\{x\in X: |g(x)h(x)|\geq\frac{1}{n}\}=\infty$ for some $n\in\mathbb{N}$. Let $u=\frac{1}{4n(1+|g|)}$. Then $u$ is a positive unit in $\mathcal{M}(X,\mathcal{A})$. We claim that $m_\mu(h,u)\cap T_g=\emptyset$ (and hence $T_g$ is closed in $\mathcal{M}_{m_\mu}(X,\mathcal{A})$). If possible let $m_\mu(h,u)\cap T_g\neq\emptyset$ and we choose $k\in m_\mu(h,u)\cap T_g$. The simple inequality $|h(x)g(x)|\leq |h(x)-k(x)||g(x)|+|k(x)||g(x)|$ for each $x\in X$ implies that: $\{x\in X:|g(x)h(x)|\geq\frac{1}{n}\}\subseteqq \{x\in X:|h(x)-k(x)|\geq u(x)\}\cup\{x\in X:|g(x)|\geq\frac{1}{2nu(x)}\}\cup\{x\in X:|k(x)g(x)|\geq \frac{1}{2n}\}\ldots(A)$. Now $u=\frac{1}{4n(1+|g|)}\implies\forall x\in X,~|g(x)|<\frac{1}{2nu(x)}\implies\{x\in X:|g(x)|\geq \frac{1}{2nu(x)}\}=\emptyset$. The inequality $(A)$ therefore implies that $\{x\in X:|g(x)h(x)|\geq\frac{1}{n}\}\subseteqq \{x\in X:|h(x)-k(x)|\geq u(x)\}\cup\{x\in X:|k(x)g(x)|\geq\frac{1}{2n}\} \ldots(B)$. Now since $k\in m_\mu(h,u)$, we get that $|k(x)-h(x)|<u(x)$ almost everywhere on $X$ and then: $\mu\{x\in X: |k(x)-h(x)|\geq u(x)\}=0$. This relation combined with $(B)$ therefore yields that $\mu\{x\in X:|g(x)h(x)|\geq\frac{1}{n}\}\leq\mu\{x\in X:|k(x)g(x)|\geq\frac{1}{2n}\}<\infty$ because $k.g\in M_\infty(X,\mathcal{A})$ as $k\in T_g$. This contradicts the relation: $\mu\{x\in X:|g(x)h(x)|\geq\frac{1}{n}\}=\infty$, already obtained earlier. Thus it is proved that $\mathcal{M}_\psi(X,\mathcal{A})$ is closed in $\mathcal{M}_{m_\mu}(X,\mathcal{A})$. To complete the Theorem let $f\in\mathcal{M}_\psi(X,\mathcal{A})$ and $u$, a positive unit in $\mathcal{M}(X,\mathcal{A})$. we need to check that $m_\mu(f,u)\cap \mathcal{M}_F(X,\mathcal{A})\neq\emptyset$. For that purpose define $g:X\to\mathbb{R}$ as follows: $$g(x)=\begin{cases}
			f(x)+\frac{\epsilon}{2}&\text{ if }f(x)\leq -\frac{\epsilon}{2}u(x)\\
			0&\text{ if }-\frac{\epsilon}{2}u(x)\leq f(x)\leq\frac{\epsilon}{2}u(x)\\
			f(x)-\frac{\epsilon}{2}&\text{ if }f(x)\geq\frac{\epsilon}{2}u(x)
		\end{cases}$$Then $g\in\mathcal{M}(X,\mathcal{A})$ and indeed $g\in m_\mu(f,u)$. We shall show that $g\in\mathcal{M}_F(X,\mathcal{A})$ and the Theorem will be completed. Indeed $X\setminus Z(g)\subseteqq \{x\in X:|f(x)|\geq\frac{1}{2}u(x)\}\implies\mu(X\setminus Z(g))\leq \mu\{x\in X:|f(x)|\geq\frac{1}{2}u(x)\}\ldots(C)$. But $f\in\mathcal{M}_\psi(X,\mathcal{A})\implies f.\frac{1}{u}\in \mathcal{M}_\infty(X,\mathcal{A})\implies\mu\{x\in X:|\frac{f(x)}{u(x)}|\geq\frac{1}{2}\}<\infty\implies\mu(X\setminus Z(g))<\infty$, from $(C)$ and hence $g\in\mathcal{M}_F(X,\mathcal{A})$.
	\end{proof}
	\begin{remark}
		Since $\mathcal{M}_F(X,\mathcal{A})$ is an ideal in $\mathcal{M}(X,\mathcal{A})$, a fact easily verifiable, and the closure of an ideal in a topological ring is very much an ideal (may be improper) [See $2M1$~\cite{Gillman7}], it follows that $\mathcal{M}_\psi(X,\mathcal{A})$ is an ideal in $\mathcal{M}(X,\mathcal{A})$ [This fact can also be established independently]. The following result highlights a sharper information about the nature of this ideal.
	\end{remark}
	\begin{theorem}\label{Th4.6}
		$\mathcal{M}_\psi(X,\mathcal{A})$ is the largest ideal in $\mathcal{M}_F(X,\mathcal{A})$ lying between $\mathcal{M}_F(X,\mathcal{A})$ and $\mathcal{M}_\infty(X,\mathcal{A})$.
	\end{theorem}
	\begin{proof}
		Let $J$ be an ideal in $\mathcal{M}(X,\mathcal{A})$ with the relation: $\mathcal{M}_F(X,\mathcal{A})\subseteqq J\subseteqq \mathcal{M}_\infty(X,\mathcal{A})$. If possible let $J\not\subset\mathcal{M}_\psi(X,\mathcal{A})$. Choose $g\in J\setminus \mathcal{M}_\psi(X,\mathcal{A})$. Then $g\notin \mathcal{M}_\psi(X,\mathcal{A})\implies g.h\notin \mathcal{M}_\infty(X,\mathcal{A})$ for some $h\in \mathcal{M}_F(X,\mathcal{A})$. But $g\in J$ and $J$ is an ideal in $\mathcal{M}_F(X,\mathcal{A})\implies g.h\in J\implies g.h\in \mathcal{M}_\infty(X,\mathcal{A})$ -- a contradiction.
	\end{proof}
	\begin{remark}
		If $U_\mu$-topology $=m_\mu$-topology, then it follows that the closures of $\mathcal{M}_F(X,\mathcal{A})$ in these two topologies are the same. Hence it follows from Theorem \ref{Th4.4} and Theorem \ref{Th4.6} that $\mathcal{M}_\psi(X,\mathcal{A})=\mathcal{M}_\infty(X,\mathcal{A})$ and consequently $\mathcal{M}_\infty(X,\mathcal{A})$ becomes an ideal in $\mathcal{M}(X,\mathcal{A})$. This implies in view of Theorem $5.3$ in \cite{Acharyya3} [where it is proved that $\mathcal{M}_\infty(X,\mathcal{A})$ is an ideal in $\mathcal{M}(X,\mathcal{A})$ if and only if $\mathcal{M}_F(X,\mathcal{A})=\mathcal{M}_\infty(X,\mathcal{A})$ that $\mathcal{M}_F(X,\mathcal{A})=\mathcal{M}_\psi(X,\mathcal{A})=\mathcal{M}_\infty(X,\mathcal{A})$].
	\end{remark}
	The following counterexample shows that, even if $U_\mu$-topology $\subsetneqq m_\mu$-topology, the relation $\mathcal{M}_F(X,\mathcal{A})=\mathcal{M}_\psi(X,\mathcal{A})=\mathcal{M}_\infty(X,\mathcal{A})$ and $\mathcal{M}_\infty(X,\mathcal{A})$ is an ideal in $\mathcal{M}(X,\mathcal{A})$ may still be valid.
	\begin{example}
		Let $X=[0,1]$, $\mathcal{A}$ be the collection of all Lebesgue measurable subsets of $[0,1]$ and $\mu=$ Lebesgue measure on it. Since $\mu(X)<\infty$, it follows form the definition that, $\mathcal{M}_F(X,\mathcal{A})=\mathcal{M}_\infty(X,\mathcal{A})$ and thus $\mathcal{M}_\infty(X,\mathcal{A})$ is an ideal in $\mathcal{M}(X,\mathcal{A})$. However, $U_\mu$-topology $\neq m_\mu$-topology as there is an $f\in \mathcal{M}(X,\mathcal{A})$, which is not essentially bounded on $X$ (vide Corollary \ref{Cor2.7}). Indeed the function $f:X\to\mathbb{R}$ defined by $f(x)=\begin{cases}
			\frac{1}{x}&\text{ if }x\neq 0\\
			0&\text{ if }x=0
		\end{cases}$ is a Lebesgue measurable function which is not essentially bounded on $X$ [The reason is for any $n\in\mathbb{N}$, $\mu(\{x\in X:f(x)\geq n\})=\mu((0,\frac{1}{n}))=\frac{1}{n}>0$].
	\end{example}
	Since Lebesgue measure on $[0,1]$ is non-atomic in the sense that each Lebesgue measurable set with positive Lebesgue measure can be decomposed into two Lebesgue measurable sets each with positive measure [$2:14.7(\text{f})$, \cite{Bruckner6}], we feel tempted to ask the following question:
	\begin{question}
		Suppose $(X,\mathcal{A},\mu)$ is a measure space with the following two conditions: \begin{enumerate}
			\item $\mu$ is an atomic measure and
			\item $\mathcal{M}_\infty(X,\mathcal{A})$ is an ideal in $\mathcal{M}(X,\mathcal{A})$.
		\end{enumerate}
		Is then $U_\mu$-topology $=m_\mu$-topology on $\mathcal{M}(X,\mathcal{A})$?
	\end{question}
	In view of Theorem \ref{Th4.0}, we would also like to make the following conjecture.\\
	\textbf{Conjecture: } The closure of $L_\psi(X,\mathcal{A},\mu)$ in the $U_\mu$-topology on $\mathcal{M}(X,\mathcal{A})$ is $L^\infty(X,\mathcal{A},\mu)$ i.e., $cl_{U_\mu}L_\psi(X,\mathcal{A},\mu)=L^\infty(X,\mathcal{A},\mu)$.
	\section{$Z_\mu$-ideals}\label{Sec5}
	We recall that an ideal $I$ in the ring $\mathcal{M}(X,\mathcal{A})$ is called a $Z_\mu$-ideal if $f\equiv 0$ almost everywhere $(\mu)$ on some zero-set $Z\in Z[I]$ implies that $f\in I$. Equivalently $I$ is a $Z_\mu$-ideal in $\mathcal{M}(X,\mathcal{A})$ if and only if $f\equiv g$ almost everywhere $(\mu)$ on $X$ and $g\in I\implies f\in I$, here $f\in \mathcal{M}(X,\mathcal{A})$.
	\begin{notation}
		Let
		\begin{align*}
			\mathcal{L}^0(X,\mathcal{A},\mu)&=\{f\in\mathcal{M}(X,\mathcal{A}): \mu(X\setminus Z(f))=0\}\\
			\mathcal{L}^F(X,\mathcal{A},\mu)&=\{f\in\mathcal{M}(X,\mathcal{A}): \mu(X\setminus Z(f))<\infty\}\\
			L_\psi(X,\mathcal{A},\mu)&=\{f\in\mathcal{M}(X,\mathcal{A}): X\setminus Z(f)\text{ is a }\mu\text{-bounded subset of }X\}
		\end{align*}
	\end{notation}
	It is easy to check that $\mathcal{L}^0(X,\mathcal{A},\mu),\mathcal{L}^F(X,\mathcal{A},\mu)$ and $L_\psi(X,\mathcal{A},\mu)$ are $Z_\mu$-ideals in this ring and also $\mathcal{L}^0(X,\mathcal{A},\mu)\subset\mathcal{L}^F(X,\mathcal{A},\mu)$, $\mathcal{L}^0(X,\mathcal{A},\mu)\subset L_\psi(X,\mathcal{A},\mu)$. More generally we have the following result.
	\begin{theorem}\label{Th5.0}
		An ideal $I$ in $\mathcal{M}(X,\mathcal{A})$ is a $Z_\mu$-ideal if and only if $\mathcal{L}^0(X,\mathcal{A},\mu)\subseteqq I$.
	\end{theorem}
	\begin{proof}
		Let $I$ be any $Z_\mu$-ideal in $\mathcal{M}(X,\mathcal{A})$ and $f\in \mathcal{L}^0(X,\mathcal{A},\mu)$. Then $f=0$ almost everywhere $(\mu)$ on $X$, because $\mu(X\setminus Z(f))=0$. Since $I$ is a $Z_\mu$-ideal and $0\in I$, this implies that $f\in I$. Thus $\mathcal{L}^0(X,\mathcal{A},\mu) \subset I$. Conversely let $I$ be an ideal in $\mathcal{M}(X,\mathcal{A})$ such that $\mathcal{L}^0(X,\mathcal{A},\mu)\subset I$. Let $f\in \mathcal{M}(X,\mathcal{A})$ be such that $f\equiv 0$ on some $Z(g)\in Z[I]$ i.e., $\mu(Z(g)\setminus Z(f))=0$. Let us consider $Z(g)\setminus Z(f)\neq\emptyset$, then there exists $h\in\mathcal{L}^F(X,\mathcal{A},\mu)$ such that $X\setminus Z(h)=Z(g)\setminus Z(f)$. Then $Z(h)=X\setminus Z(g)\cup Z(f)\implies Z(f)\cap Z(g)=Z(g)\cap Z(h)$ i.e., $Z(f^2+g^2)=Z(g^2+h^2)\implies f^2=(g^2+h^2).f_1-g^2$ for some $f_1\in \mathcal{M}(X,\mathcal{A})$ [Theorem $2.3$, \cite{Acharyya1}]. Since $\mathcal{L}^0(X,\mathcal{A},\mu)\subset I$ and $h\in \mathcal{L}^0(X,\mathcal{A},\mu)$, then $h\in I\implies g^2+h^2\in I$. Thus $f^2\in I\implies f\in I$. Consequently, $I$ is a $Z_\mu$-ideal.
	\end{proof}
	\begin{corollary}\label{Rem5.2}
		$\mathcal{L}^0(X,\mathcal{A},\mu)$ is the smallest $Z_\mu$-ideal in $\mathcal{M}(X,\mathcal{A})$.
	\end{corollary}
	\begin{corollary}\label{Th5.199}
		Let $M$ be a maximal ideal in $\mathcal{M}(X,\mathcal{A})$ containing a $Z_\mu$-ideal $I$. Then $M$ is a $Z_\mu$-ideal. In particular, a maximal $Z_\mu$-ideal in $\mathcal{M}(X,\mathcal{A})$ is a maximal ideal.
	\end{corollary}
	The intersection of any collection of $Z_\mu$-ideals in $\mathcal{M}(X,\mathcal{A})$ is evidently a $Z_\mu$-ideal. Hence for any ideal $I$ in $\mathcal{M}(X,\mathcal{A})$, $I^e\equiv$ the intersection of all the $Z_\mu$-ideals which contain $I$, is the smallest $Z_\mu$-ideal containing $I$. The following result offers a convenient formulae for $I^e$.
	\begin{theorem}\label{Th5.1}
		Let $I$ be an ideal in $\mathcal{M}(X,\mathcal{A})$. If there exists an $f\in I$ such that $\mu(Z(f))=0$, then $I^e=\mathcal{M}(X,\mathcal{A})$. If however for each $f\in I$, $\mu(Z(f))>0$, then $I^e=\{g\in\mathcal{M}(X,\mathcal{A}):\exists~f\in I\text{ such that }\mu(Z(f)\triangle Z(g))=0\}$. [Here $A\triangle B$ denotes the symmetric difference between two sets $A,B$]
	\end{theorem}
	\begin{proof}
		First let there exist $f\in I$ such that $\mu(Z(f))=0$. Choose $g\in\mathcal{M}(X,\mathcal{A})$ and an arbitrary $Z_\mu$-ideal $J$ containing $I$. Then $Z(f)\in Z[I]$ and $\mu(Z(f))=0$ implies that $g\in J$, because the condition $g\equiv 0$ almost everywhere $(\mu)$ on  $Z(f)\in Z[I]$ is trivially satisfied. Hence $g\in I^e$. Thus $I^e=\mathcal{M}(X,\mathcal{A})$. Now let $\forall f\in I$, $\mu(Z(f))>0$. Suppose $\mathcal{B}=\{A\in\mathcal{A}:\exists~f\in I\text{ such that }\mu(Z(f)\triangle A)=0\}$. Since all the measurable sets in the $\mathcal{A}$-filter (a filter of $\mathcal{A}$-measurable sets) $Z[I]=\{Z(f):f\in I\}$ are thick measurable sets in the sense that for each $Z\in Z[I]$, $\mu(Z)>0$, it is easy to easy that $\emptyset\notin\mathcal{B}$. Furthermore it is not at all hard to check that if $B_1,B_2\in\mathcal{B}$ and $B_1\subseteq B_3\in\mathcal{A}$, then $B_1\cap B_2\in\mathcal{B}$ and $B_3\in\mathcal{B}$. Thus $\mathcal{B}$ is an $\mathcal{A}$-filter on $X$ and surely $Z[I]\subseteqq\mathcal{B}$. Let $Z^{-1}[\mathcal{B}]=J$ i.e., $Z[J]=\mathcal{B}$. Then $J$ is an ideal in $\mathcal{M}(X,\mathcal{A})$ containing $I$. To complete the Theorem it suffices to show that $J=I^e$ i.e., $J$ is the intersection of all the $Z_\mu$-ideals in $\mathcal{M}(X,\mathcal{A})$ containing $I$. First we show that $J$ is a $Z_\mu$-ideal. Let $f\in \mathcal{M}(X,\mathcal{A})$ such that $f\equiv g$ almost everywhere $(\mu)$ for some $g\in J$ i.e., $\mu(Z(g)\setminus Z(f))=0$ and $\mu(Z(f)\setminus Z(g))=0\implies \mu(Z(f)\triangle Z(g))=0$. Since $g\in J$, there exists $h\in I$ such that $\mu(Z(g)\triangle Z(h))=0$. Now $Z(f)\triangle Z(h)\subset(Z(f)\triangle Z(g))\cup (Z(g)\triangle Z(h))\implies \mu(Z(f)\triangle Z(h))=0\implies Z(f)\in\mathcal{B}$ i.e., $f\in J$. Hence $J$ is a $Z_\mu$-ideal. Let $K$ be any $Z_\mu$-ideal in $\mathcal{M}(X,\mathcal{A})$ containing $I$. We assert that $J\subseteqq K$. Proof of the assertion: let $f\in J$. Then there exists $g\in I$ such that $\mu(Z(f)\triangle Z(g))=0\implies \mu(Z(g)\setminus Z(f))=0$ i.e., $f\equiv 0$ almost everywhere $(\mu)$ on $Z(g)\in Z[I]$. Since $K$ is a $Z_\mu$-ideal containing $I$, $g\in K$ and $f\equiv 0$ almost everywhere $(\mu)$ on $Z(g)\in Z[K]\implies f\in K$. Hence $J$ is the intersection of all the $Z_\mu$-ideals in $\mathcal{M}(X,\mathcal{A})$ containing $I$ and that completes the proof.
	\end{proof}
	A maximal ideal in $\mathcal{M}(X,\mathcal{A})$ may or may not be a $Z_\mu$-ideal in $\mathcal{M}(X,\mathcal{A})$. Consider the following examples:
	\begin{example}\label{Ex5.3}
		\textit{(Example of Maximal ideals which are also $Z_\mu$-ideals)}
		\begin{enumerate}
			\item Let $\mathscr{F}=\{A\in\mathcal{A}:\mu(X\setminus A)=0\}$. It is easy to check that $\mathcal{F}$ is a filter of $\mathcal{A}$-measurable sets in $X$. Hence $\mathscr{F}$ can be extended to an $\mathcal{A}$-ultrafilter $\mathscr{U}$ of $\mathcal{A}$-measurable sets in $X$. Consequently by Theorem $2.7$ in \cite{Acharyya1}, $\mathscr{U}=Z[M]\equiv\{Z(f):f\in M\}$ for some maximal ideal $M$ in $\mathcal{M}(X,\mathcal{A})$. We assert that for each $A\in\mathscr{U}=Z[M]$, $\mu(A)>0$. Proof of this assertion: if possible let for some $A\in Z[M]$, $\mu(A)=0$. Then $X\setminus A\in \mathscr{F}\subset\mathscr{U}$. Therefore $\emptyset=X\setminus A\cap A\in\mathscr{U}$ -- a contradiction. We now show that $M^e\equiv$ the smallest $Z_\mu$-ideal in $\mathcal{M}(X,\mathcal{A})$ containing $M$ is a proper ideal in $\mathcal{M}(X,\mathcal{A})$. If not then $1\in M^e$ and hence by Theorem \ref{Th5.1} there exists $f\in M$ such that $\mu(Z(f)\triangle Z(1))=0\implies \mu(Z(f))=0$ - a contradiction to the fact that $\mu(Z(f))>0$. Since $M$ is a maximal ideal in $\mathcal{M}(X,\mathcal{A})$, $M=M^e$ i.e., $M$ is a $Z_\mu$-ideal.
			\item Consider the $\mathcal{A}$-filter $\mathscr{F'}=\{A\in\mathcal{A}: \mu(X\setminus A)<\infty\}$. Similarly, $\mathscr{F'}$ can be extended to an $\mathcal{A}$-ultrafilter $\mathscr{U'}$ where $\mu(A)=\infty$ for all $A\in\mathscr{U'}$. Let $M'$ be a maximal ideal in $\mathcal{M}(X,\mathcal{A})$ such that $Z[M']=\mathscr{U'}$. Following the same approach as above, we can show that $M'$ is a $Z_\mu$-ideal.
		\end{enumerate}
	\end{example}
	\begin{example}\label{Ex5.4}
		\textit{(Example of a Maximal ideal which is not a $Z_\mu$-ideal)}\\
			Suppose $(X,\mathcal{A},\mu)$ be a measure space such that there exists $x\in X$ with $\{x\}\in\mathcal{A}$ and $\mu(\{x\})=0$ (For example, $(X,\mathcal{A},\mu)\equiv$ Lebesgue measure on $\mathbb{R}$). Consider the fixed maximal ideal $M_x=\{f\in \mathcal{M}(X,\mathcal{A}):f(x)=0\}$ in $\mathcal{M}(X,\mathcal{A})$. Then $\{x\}\in Z[M_x]$ and $\mu(\{x\})=0$. By Theorem \ref{Th5.1}, $M_x^e=\mathcal{M}(X,\mathcal{A})$ and hence $M_x$ is not a $Z_\mu$-ideal.
	\end{example}
	We recall from \cite{Acharyya2}, that an $u\in \mathcal{M}(X,\mathcal{A})$ is called a $\mu$-unit if $\mu(Z(u))=0$. It is also proved in Theorem $3.2$ in \cite{Acharyya2} that the set of all $\mu$-units in $\mathcal{M}(X,\mathcal{A})$ is open in the $m_\mu$-topology. Also is is easy to prove that no $\mu$-unit in $\mathcal{M}(X,\mathcal{A})$ can be a member of a proper $Z_\mu$-ideal. Hence the $m_\mu$-closure of a proper $Z_\mu$-ideal in $\mathcal{M}(X,\mathcal{A})$ is also a proper ideal. This fact combined with the next result implies that each maximal $Z_\mu$-ideal in $\mathcal{M}(X,\mathcal{A})$ is closed in $m_\mu$-topology. 
	\begin{theorem}\label{Th5.2}
		For any ideal $I$ in $\mathcal{M}(X,\mathcal{A})$, its closure in the $m_\mu$-topology is a $Z_\mu$-ideal.
	\end{theorem}
	\begin{proof}
		Since $\mathcal{M}(X,\mathcal{A})$ equipped with the $m_\mu$-topology is a topological ring, the closure of $I$ in the $m_\mu$-topology i.e., $cl_{m_\mu}I$ is an ideal in $\mathcal{M}(X,\mathcal{A})$. Let $g\in cl_{m_\mu}I$ and $f\equiv g$ almost everywhere $(\mu)$ on $X$, where $f\in \mathcal{M}(X,\mathcal{A})$. We shall show that $f\in cl_{m_\mu}I$. Let $u$ be a positive unit in $\mathcal{M}(X,\mathcal{A})$. Then $m_\mu(g,u)\cap I\neq\emptyset$, as $g\in cl_{m_\mu}I$. Choose $h\in m_\mu(g,u)\cap I$. Then there exists $A\in\mathcal{A}$ with $\mu(A)=0$ such that $|h-g|<u$ on $X\setminus A$. Suppose $Z=\{x\in X:f(x)\neq g(x)\}$, then $\mu(Z)=0$. It follows that $\mu(A\cup Z)=0$ and $|h-f|=|h-g|<u$ on $X\setminus (A\cup Z)$. Therefore $h\in m_\mu(f,u)\cap I$. Hence $f\in cl_{m_\mu}I$.
	\end{proof}
	\begin{corollary}\label{Cor5.3}
		For any ideal $I$ in $\mathcal{M}(X,\mathcal{A})$, $I\subseteqq I^e \subseteqq cl_{m_\mu}I$.
	\end{corollary}
	\begin{corollary}\label{Cor5.8}
		A closed ideal in $\mathcal{M}(X,\mathcal{A})$ in the $m_\mu$-topology is necessarily a $Z_\mu$-ideal
	\end{corollary}
	Since $\mathcal{M}(X,\mathcal{A})$ equipped with the $U_\mu$-topology is not a topological ring, the closure of an ideal in the $U_\mu$-topology may not be an ideal in $\mathcal{M}(X,\mathcal{A})$. But if $I$ is such an ideal in $\mathcal{M}(X,\mathcal{A})$ whose $U_\mu$-closure is also an ideal in $\mathcal{M}(X,\mathcal{A})$, then the closure of $I$ in the $U_\mu$-topology i.e., $cl_{U_\mu}I$ is also a $Z_\mu$-ideal. We can easily establish this result by closely following the proof of Theorem \ref{Th5.2}. Nevertheless since the $m_\mu$-topology on $\mathcal{M}(X,\mathcal{A})$ is finer than the $U_\mu$-topology, we have the following result:
	\begin{remark}\label{Rem5.6}
		For any ideal $I$ in $\mathcal{M}(X,\mathcal{A})$, $I\subseteqq I^e \subseteqq cl_{m_\mu}I\subseteqq cl_{U_\mu}I$.
	\end{remark}
	The following Theorem characterizes the proper ideals in $\mathcal{M}(X,\mathcal{A})$ which are dense in the $m_\mu$-topology/ $U_\mu$-topology.
	\begin{theorem}\label{Th5.6}
		For a proper ideal $I$ in $\mathcal{M}(X,\mathcal{A})$, the following statements are equivalent:
		\begin{enumerate}
			\item There exists $Z\in Z[I]$ such that $\mu(Z)=0$.
			\item $I$ is dense in $\mathcal{M}(X,\mathcal{A})$ in the $m_\mu$-topology.
			\item $I$ is dense in $\mathcal{M}(X,\mathcal{A})$ in the $U_\mu$-topology.
		\end{enumerate}
	\end{theorem}
	\begin{proof}
		Let $(1)$ hold i.e., $\mu(Z)=0$ for some $Z\in Z[I]$. Then by Theorem \ref{Th5.1}, $I^e=\mathcal{M}(X,\mathcal{A})$ and hence by Corollary \ref{Cor5.3}, $cl_{m_\mu}I=\mathcal{M}(X,\mathcal{A})$. Therefore $(1)\implies(2)$. $(2)\implies (3)$ is trivial because the $m_\mu$-topology on $\mathcal{M}(X,\mathcal{A})$ is finer than the $U_\mu$-topology. Suppose $(3)$ is true i.e., $cl_{U_\mu}I=\mathcal{M}(X,\mathcal{A})$. So $1\in cl_{U_\mu}I\implies U_\mu(1,\frac{1}{2})\cap I\neq\emptyset$. Choose $g\in U_\mu(1,\frac{1}{2})\cap I$ i.e., there exists $A\in\mathcal{A}$ with $\mu(A)=0$ such that $|g-1|<\frac{1}{2}$ on $X\setminus A$. This means that $\frac{1}{2}<g(x)<\frac{3}{2}$ for all $x\in X\setminus A\implies Z(g)\subset A\implies \mu(Z(g))=0$ and also $g\in I$. Thus $(3)\implies (1)$
	\end{proof}
	\begin{corollary}
		A proper ideal $I$ in $\mathcal{M}(X,\mathcal{A})$ is dense in the $m_\mu$-topology/ $U_\mu$-topology if and only if there exists a $\mu$-unit in $I$.
	\end{corollary}
	\begin{corollary}
		A proper $Z_\mu$-ideal in $\mathcal{M}(X,\mathcal{A})$ in the $m_\mu$-topology/ $U_\mu$-topology is never dense in either spaces. 
	\end{corollary}
	It follows from Remark \ref{Rem5.6} that if an ideal in $\mathcal{M}(X,\mathcal{A})$ is closed in the $m_\mu$-topology, then the ideal is necessarily a $Z_\mu$-ideal. We now show that the converse is also true i.e., every $Z_\mu$-ideal in $\mathcal{M}(X,\mathcal{A})$ is closed in the $m_\mu$-topology. Hence the set of all $Z_\mu$-ideals in $\mathcal{M}(X,\mathcal{A})$ is precisely the set of all closed ideal in the $m_\mu$-topology.  
	\begin{notation}
		For any ideal $I$ in $\mathcal{M}(X,\mathcal{A})$, let $I_\mu=\{f\in\mathcal{M}(X,\mathcal{A}):\forall~g\in\mathcal{M}(X,\mathcal{A}),~f.g\in cl_{U_\mu}I\}$.
	\end{notation}
	We establish the following facts successively.
	\begin{theorem}\label{Th5.13}
		Let $I$ be an ideal in $\mathcal{M}(X,\mathcal{A})$. Then,
		\begin{enumerate}
			\item $I_\mu$ is an ideal in $\mathcal{M}(X,\mathcal{A})$ containing $I$.
			\item $I_\mu$ is the largest ideal containing $I$ and contained in $cl_{U_\mu}I$.\label{2}
			\item $I_\mu$ is a closed subset of $\mathcal{M}(X,\mathcal{A})$ in the $m_\mu$-topology.\label{3}
			\item $I_\mu=cl_{m_\mu}I$.
			\item $I_\mu$ is a $Z_\mu$-ideal.\label{5}
			\item $I^e=I_\mu$.\label{6}
			\item $I$ is a $Z_\mu$-ideal if and only if $I=I_\mu$.
		\end{enumerate}
	\end{theorem}
	\begin{proof}
		\hspace*{3cm}
		\begin{enumerate}
			\item Since $cl_{U_\mu}I$ is a additive subgroup of $\mathcal{M}(X,\mathcal{A})$, it is easy to show that $I_\mu$ is an ideal in $\mathcal{M}(X,\mathcal{A})$ and $I\subseteqq I_\mu$.
			\item Let $J$ be an ideal in $\mathcal{M}(X,\mathcal{A})$ such that $I\subset J\subset cl_{U_\mu}I$ and $f\in J$. Then for any $g\in \mathcal{M}(X,\mathcal{A})$, $f.g\in J\subset cl_{U_\mu}I\implies f\in I_\mu$ i.e., $J\subset I_\mu$. Thus $I_\mu$ is the largest ideal lying between $I$ and $cl_{U_\mu}I$.
			\item Let $f\in \mathcal{M}(X,\mathcal{A})\setminus I_\mu$. Then there exist $g\in \mathcal{M}(X,\mathcal{A})$ such that $f.g\notin cl_{U_\mu}I\implies$ there exists $\epsilon>0$ such that $U_\mu(f.g,\epsilon)\cap I=\emptyset$. Let $u=\frac{\epsilon}{2(1+|g|)}$, then we assert that $m_\mu(f,u)\cap I_\mu=\emptyset$. If possible let there exists $h\in m_\mu(f,u)\cap I_\mu$. Then $|h-f|<u$ almost everywhere $(\mu)$ on $X\implies |h.g-f.g|<u|g|<\epsilon$ almost everywhere $(\mu)$ on $X\implies h.g\in U_\mu(f.g,\epsilon)$. Now $h\in I_\mu\implies h.g\in cl_{U_\mu}I\implies U_\mu(f.g,\epsilon)\cap I\neq\emptyset$ - contradiction.
			\item Since $cl_{m_\mu}I$ is the smallest closed set in $\mathcal{M}(X,\mathcal{A})$ containing $I$ and $I_\mu(\supseteqq I)$ is a closed set in $\mathcal{M}(X,\mathcal{A})$ in the $m_\mu$-topology (from Theorem \ref{Th5.13}(\ref{3})), it follows that $cl_{m_\mu}I\subseteqq I_\mu$. To prove the reverse implication relation let $f\in I_\mu$ and $u$, a positive unit in $\mathcal{M}(X,\mathcal{A})$. To show that $m_\mu(f,u)\cap I\neq\emptyset$, take $g=\frac{2}{u}$. Then $f\in I_\mu\implies f.g\in cl_{U_\mu}I\implies U_\mu(f.g,1)\cap I\neq\emptyset$. Choose $h\in U_\mu(f.g,1)\cap I$. Then $h\in I$ and $|h-f.g|<1$ almost everywhere $(\mu)$ on $X\implies |h.u-2f|=|h.u-f.g.u|<u$ almost everywhere $(\mu)$ on $X\implies |\frac{h.u}{2}-f|<\frac{u}{2}<u$ almost everywhere $(\mu)$ on $X\implies \frac{h.u}{2}\in m_\mu(f,u)\implies\frac{h.u}{2}\in m_\mu(f,u)\cap I$.
			\item We already have $I_\mu=cl_{m_\mu}I$ and hence from Theorem \ref{Th5.2} it follows that $I_\mu$ is a $Z_\mu$-ideal.
			\item Since $I^e$ is the smallest $Z_\mu$-ideal containing $I$, it follows from Theorem \ref{Th5.13}(\ref{5}) that $I^e\subseteqq I_\mu$. Let $f\in I_\mu$. Consider $g(x)=\begin{cases}
				\frac{1}{f(x)}&\text{ if }x\notin Z(f)\\
				0&\text{ if }x\in Z(f)
			\end{cases}$. Then $f.g\in cl_{U_\mu}I\implies U_\mu(f.g,\frac{1}{2})\cap I\neq\emptyset$. Suppose $h\in U_\mu(f.g,\frac{1}{2})\cap I$. Then there exists $A\in\mathcal{A}$ with $\mu(A)=0$ such that $|h-f.g|<\frac{1}{2}$ on $X\setminus A$. Let $x\in Z(h)\cap X\setminus A$. Then $h(x)=0$ and $|h(x)-f(x).g(x)|<\frac{1}{2}\implies |f(x).g(x)|<\frac{1}{2}\implies x\in Z(f)$, otherwise $f(x).g(x)=1$. Therefore $Z(h)\cap (X\setminus A)\subseteqq Z(f)\implies f\equiv 0$ almost everywhere $(\mu)$ on $Z(h)\in Z[I]\subseteqq Z[I^e]\implies f\in I^e$. Hence $I^e=I_\mu$.
			\item Immediate consequence of Theorem \ref{Th5.13}(\ref{6}) and the fact that $I$ is a $Z_\mu$-ideal if and only if $I=I^e$.
		\end{enumerate}
	\end{proof}
	\begin{remark}
		For any ideal $I$ in $\mathcal{M}(X,\mathcal{A})$, $cl_{m_\mu}I= \{f\in\mathcal{M}(X,\mathcal{A}):\forall~g\in\mathcal{M}(X,\mathcal{A}),~f.g\in cl_{U_\mu}I\}$ and $cl_{m_\mu}I$ is the largest ideal lying between $I$ and $cl_{U_\mu}I$.
	\end{remark}
	\begin{corollary}
		For any ideal $I$ in $\mathcal{M}(X,\mathcal{A})$, $cl_{U_\mu}I$ is an ideal in $\mathcal{M}(X,\mathcal{A})$ if and only if $cl_{m_\mu}I=cl_{U_\mu}I$.
	\end{corollary}
	\begin{corollary}\label{Cor5.16}
		An ideal $I$ in $\mathcal{M}(X,\mathcal{A})$ is closed in the $m_\mu$-topology if and only if it is a $Z_\mu$-ideal.
	\end{corollary}
	Hence form Theorem \ref{Th5.0}, the following result holds.
	\begin{theorem}
		$\{I\subset\mathcal{M}(X,\mathcal{A}):\mathcal{L}^0(X,\mathcal{A},\mu) \subseteqq I\text{ and }I\text{ is an ideal in }\mathcal{M}(X,\mathcal{A})\}$ is precisely the set of all ideals in $\mathcal{M}(X,\mathcal{A})$ which are closed in the $m_\mu$-topology.
	\end{theorem}
	\begin{theorem}\label{Th5.19}
		A maximal ideal $M$ in $\mathcal{M}(X,\mathcal{A})$ is a $Z_\mu$-ideal if and only if $\mu(Z(f))>0$ for all $f\in M$.
	\end{theorem}
	\begin{proof}
		If a maximal ideal $M$ is a $Z_\mu$-ideal, then it is closed in the $m_\mu$-topology (by Corollary \ref{Cor5.16}) i.e., $M$ is not Dense in the $m_\mu$-topology. Hence by Theorem \ref{Th5.6}, $\mu(Z(f))>0$ for all $f\in M$. Conversely if $\mu(Z(f))>0$ for all $f\in M$, then $M$ is not dense in the $m_\mu$-topology i.e., $cl_{m_\mu}M\subsetneqq \mathcal{M}(X,\mathcal{A})$. Since $\mathcal{M}(X,\mathcal{A})$ in the $\mu$-topology is a topological ring [$2M1$, \cite{Gillman7}], it follows that $cl_{m_\mu}M$ is a proper ideal in $\mathcal{M}(X,\mathcal{A})$. The maximality of $M$ forces that $M=cl_{m_\mu}M$ i.e., $M$ is a $Z_\mu$-ideal in $\mathcal{M}(X,\mathcal{A})$ (by Theorem \ref{Th5.2}).
	\end{proof}
	\begin{corollary}\label{Cor5.20}
		A maximal ideal $M$ in $\mathcal{M}(X,\mathcal{A})$ is closed in the $m_\mu$-topology if and only if $\mu(Z(f))>0$ for all $f\in M$.
	\end{corollary}
	\begin{theorem}\label{Th5.20}
		Let $(X,\mathcal{A},\mu)$ be a measure space where $\{x\}\in\mathcal{A}$ for all $x\in X$. Then the following two statements are equivalent:
		\begin{enumerate}
			\item Each $Z_\mu$-ideal in $\mathcal{M}(X,\mathcal{A})$ is free.
			\item $\mu(\{x\})=0$ for all $x\in X$. 
		\end{enumerate}
	\end{theorem}
	\begin{proof}
		\underline{$(1)\implies (2): $} If possible let there exist $x\in X$ such that $\mu(\{x\})>0$. We assert that $M_x=\{f\in \mathcal{M}(X,\mathcal{A}):f(x)=0\}$, which is evidently a fixed maximal ideal in $\mathcal{M}(X,\mathcal{A})$ is a $Z_\mu$-ideal. For all $f\in M_x$, $x\in Z(f)\implies \mu(Z(f))>0$. Then by Theorem \ref{Th5.19}, $M_x$ is a $Z_\mu$-ideal -- contradicts $(1)$.\\
		\underline{$(2)\implies(1): $} Let $I$ be a $Z_\mu$-ideal in $\mathcal{M}(X,\mathcal{A})$ and $x\in X$. Then $1_{\{x\}}\equiv 0$ almost everywhere $(\mu)$ on $X=Z(0)\in Z[I]\implies 1_{\{x\}}\in I$. This proves that $I$ is free.
	\end{proof}
	\section{Description of the $U^I_{\mu,F}$-topology and the $m^I_{\mu,F}$-topology on $\mathcal{M}(X,\mathcal{A})$}\label{Sec6}
	Throughout the section $I$ always stands for an ideal in $\mathcal{M}(X,\mathcal{A})$. We first state a simple comparison between all the topologies mentioned on $\mathcal{M}(X,\mathcal{A})$.
	\begin{theorem}
		On $\mathcal{M}(X,\mathcal{A})$, $U^I_\mu$-topology $\supseteqq U^I_{\mu,F}$-topology and $m^I_\mu$-topology $\supseteqq m^I_{\mu,F}$-topology on $\mathcal{M}(X,\mathcal{A})$. Also $U^I_{\mu,F}$-topology $\subseteqq m^I_{\mu,F}$-topology on $\mathcal{M}(X,\mathcal{A})$.
	\end{theorem}
	\begin{example}
		In the co-countable measure space $(X,\mathcal{A},\rho)$ (discussed in \ref{CoCu}), $U^I_\rho$-topology $=U^I_{\rho,F}$-topology and $m^I_\rho$-topology $=m^I_{\rho,F}$-topology on $\mathcal{M}(X,\mathcal{A})$ for every ideal $I$ in $\mathcal{M}(X,\mathcal{A})$ [see Theorem \ref{Th7.5}]. 
	\end{example}
	We let $\mathcal{M}_{U^I_{\mu,F}}(X,\mathcal{A})$ to stand for $\mathcal{M}(X,\mathcal{A})$ equipped with the $U^I_{\mu,F}$-topology with an analogous meaning for $\mathcal{M}_{m^I_{\mu,F}} (X,\mathcal{A})$. It can be realized without any difficulty that, most of the results related to the $U^I_{\mu,F}$-topology and the $m^I_{\mu,F}$-topology on $\mathcal{M}(X,\mathcal{A})$ have their counterparts in the $U^I_{\mu}$-topology and the $m^I_{\mu}$-topology in the previous sections. Since the proof of those results are also parallel to the proofs of their corresponding counterparts, we simply omit them. However we state all these parallel facts for our convenience.
	\begin{theorem}
		Let $I,J$ be two distinct ideals in $\mathcal{M}(X,\mathcal{A})$. Then on $\mathcal{M}(X,\mathcal{A})$, $U^I_{\mu,F}$-topology $\neq U^J_{\mu,F}$-topology and $m^I_{\mu,F}$-topology $\neq m^J_{\mu,F}$-topology.
	\end{theorem}
	\begin{theorem}
		If $I$ and $J$ are two ideals in $\mathcal{M}(X,\mathcal{A})$ such that $I\subseteqq J$, then the $U^I_{\mu,F}$-topology is finer than the $U^J_{\mu,F}$-topology on $\mathcal{M}(X,\mathcal{A})$ (and also $m^J_{\mu,F}$-topology $\subseteqq m^I_{\mu,F}$-topology on $\mathcal{M}(X,\mathcal{A})$).
	\end{theorem}
	\begin{theorem}
		The following statements are equivalent for an ideal $I$ in $\mathcal{M}(X,\mathcal{A})$:
		\begin{enumerate}
			\item $I=\{0\}$.
			\item $\mathcal{M}_{U^I_{\mu,F}}(X,\mathcal{A})$ is a discrete space.
			\item $\mathcal{M}_{m^I_{\mu,F}}(X,\mathcal{A})$ is a discrete space.
		\end{enumerate}
	\end{theorem}
	\begin{definition}
		An $f\in\mathcal{M}(X,\mathcal{A})$ is called weakly essentially bounded in $X$ if $f$ is bounded except for a $\mu$-finite set i.e., there exists $A_f\in\mathcal{A}$ such that $|f|<\lambda$ on $X\setminus A_f$ for some $\lambda>0$.
	\end{definition}
	\begin{notation}
		$L^\infty_F(X,\mathcal{A},\mu)=\{f\in\mathcal{M}(X,\mathcal{A}):f\text{ is weakly essentially bounded in }X\}$
	\end{notation}
	It is easy to show that $L^\infty_F(X,\mathcal{A},\mu)$ is a subring of $\mathcal{M}(X,\mathcal{A})$ containing $L^\infty(X,\mathcal{A},\mu)$.
	\begin{theorem}\label{Th6.6}
		If $\mathcal{M}_{U^I_{\mu,F}}(X,\mathcal{A})$ is a topological ring, then $I\subseteqq L^\infty_F(X,\mathcal{A},\mu)$.
	\end{theorem}
	\begin{proof}
		If possible let there exist an $f\in I\setminus L^\infty_F(X,\mathcal{A},\mu)$. Without loss of generality let $f\geq 0$ on $X$. Since $\mathcal{M}_{U^I_{\mu,F}}(X,\mathcal{A})$ is a topological ring, there exists $\epsilon>0$ such that $U_{\mu,F}(0,I,\epsilon)\times U_{\mu,F}(f,I,\epsilon)\subseteqq U_{\mu,F}(0,I,1)$. Let $g=\frac{\epsilon.f}{2(1+f)}$. Since $f\in I$, $g\in I$. Also $g\geq 0$ and $|g|<\frac{\epsilon}{2}$ on $X\implies g\in U_{\mu,F}(0,I,\epsilon)\implies g.f\in U_{\mu,F}(0,I,1)$. Then there exists $A\in\mathcal{A}$ with $\mu(A)<\infty$ such that $|f.g|<1$ on $X\setminus A$ i.e., for all $x\in X\setminus A$, $\frac{\epsilon f^2(x)}{2(1+f(x))}<1$. Since $f\notin L^\infty_F(X,\mathcal{A},\mu)$, for each $n\in\mathbb{N}$, $\mu(\{x\in X\setminus A: f(x)\geq n\})=\infty$. Choose $x_n\in \{x\in X\setminus A:f(x)\geq n\}$ for each $n\in\mathbb{N}$. Then $f(x_n)\geq n$ for all $n\in\mathbb{N}\implies \lim\limits_{n\to\infty}f(x_n)=\infty\implies \lim\limits_{n\to\infty}\frac{f(x_n)}{1+f(x_n)}=1$ i.e., $\lim\limits_{n\to\infty}g(x_n)= \frac{\epsilon}{2}$. Then there exists $p\in\mathbb{N}$ such that $g(x_n)>\frac{\epsilon}{4}$ for all $n\geq p\implies f(x_n).g(x_n)\geq \frac{\epsilon}{4}f(x_n)$. Again since $\{x_n:n\geq p\}\subset X\setminus A$, $f(x_n).g(x_n)= |f(x_n).g(x_n)|<1$ for all $n\geq p\implies f(x_n)<\frac{4}{\epsilon}$ for all $n\geq p$ -- this contradicts the fact that $\lim\limits_{n\to\infty}f(x_n)=\infty$.
	\end{proof}
	\begin{corollary}
		If $U^I_{\mu,F}$-topology $=m^I_{\mu,F}$-topology on $\mathcal{M}(X,\mathcal{A})$, then $I\subseteqq L^\infty_F(X,\mathcal{A},\mu)$.
	\end{corollary}
	To find out when $\mathcal{M}(X,\mathcal{A})$ is a topological vector space in the topology $m^I_{\mu,F}$ (and $U^I_{\mu,F})$, we first prove the following lemmas.
	\begin{lemma}\label{Lem6.8}
		If $\mathcal{M}(X,\mathcal{A})=L^\infty_F(X,\mathcal{A},\mu)$, then $U^I_{\mu,F}$-topology $=m^I_{\mu,F}$-topology on $\mathcal{M}(X,\mathcal{A})$, for any ideal $I$ in $\mathcal{M}(X,\mathcal{A})$.
	\end{lemma}
	\begin{proof}
		Let $\mathcal{M}(X,\mathcal{A})=L^\infty_F(X,\mathcal{A},\mu)$ and consider any basic open set $m_{\mu,F}(f,I,u)$ in the $m^I_{\mu,F}$-topology on $\mathcal{M}(X,\mathcal{A})$, where $f\in \mathcal{M}(X,\mathcal{A})$ and $u$ is a positive unit in $\mathcal{M}(X,\mathcal{A})$. We have to show that $m_{\mu,F}(f,I,u)$ is open in the $U^I_{\mu,F}$-topology. Since $\frac{1}{u}\in\mathcal{M}(X,\mathcal{A})=L^\infty_F(X,\mathcal{A},\mu)$, there exists a $\mu$-finite set $A$ such that $|\frac{1}{u}|<\lambda$ on $X\setminus A$ for some $\lambda>0$ i.e., for all $x\in X\setminus A$, $u(x)>\frac{1}{\lambda}$. It is easy to check that $U_{\mu,F}(f,I,\frac{1}{\lambda})\subset m_{\mu,F}(f,I,u)$ and hence $m_{\mu,F}(f,I,u)$ is an open set in the $U^I_{\mu,F}$-topology. 
	\end{proof}
	\begin{lemma}\label{Lem6.9}
		\hspace*{3cm}
		\begin{enumerate}
			\item $\mathcal{M}_{U_{\mu,F}}(X,\mathcal{A})$ is a topological vector space if and only if $\mathcal{M}(X,\mathcal{A})=L^\infty_F(X,\mathcal{A},\mu)$.
			\item $\mathcal{M}_{m_{\mu,F}}(X,\mathcal{A})$ is a topological vector space if and only if $\mathcal{M}(X,\mathcal{A})=L^\infty_F(X,\mathcal{A},\mu)$.
		\end{enumerate}
	\end{lemma}
	\begin{proof}
		\hspace*{3cm}
		\begin{enumerate}
			\item Let $\mathcal{M}(X,\mathcal{A})=L^\infty_F(X,\mathcal{A},\mu)$. We have to show that the scalar multiplication \begin{alignat*}{1}
				S:\mathbb{R}\times\mathcal{M}_{U_{\mu,F}}(X,\mathcal{A})&\to \mathcal{M}_{U_{\mu,F}}(X,\mathcal{A})\\
				(r,f)&\mapsto r.f
			\end{alignat*} is continuous. To prove that the scalar multiplication is continuous at $(r,f)$, consider any basic open set $U_{\mu,F}(r.f, \epsilon)$ at $S(r,f)$. We have to find a suitable $\delta>0$ such that $S((r-\delta,r+\delta), U_{\mu,F}(f,\delta))\subseteqq U_{\mu,F}(r.f, \epsilon)$. Since $f\in L^\infty_F(X,\mathcal{A},\mu)$, there exists $\mu$-finite set $A_f$ such that $|f|<\lambda$ on $X\setminus A_f$ for some $\lambda>0$. Let $s\in (r-\delta,r+\delta)$ and $g\in U_{\mu,F}(f,\delta))$. Then $|s.g-r.f|\leq |s||g-f|+|s-r||f|\leq (|r|+\delta)\delta+\delta.\lambda=\delta(|r|+\delta+\lambda)\leq \delta(|r|+\lambda+1)$ (suppose $\delta<1$) on $X\setminus A_f$. If we choose $\delta=\frac{\epsilon}{2(|r|+\lambda+1)}\wedge 1$, then $\delta<1$ and $\delta(|r|+\lambda+1)\leq\frac{\epsilon}{2} <\epsilon$ i.e., $|s.g-r.f|<\epsilon$ on $X\setminus A_f\implies s.g\in U_{\mu,F}(r.f,\epsilon)$. Consequently, $\mathcal{M}_{m_{\mu,F}}(X,\mathcal{A})$ is a topological vector space. Conversely let $\mathcal{M}_{U_{\mu,F}}(X,\mathcal{A})$ be a topological vector space and $f\in\mathcal{M}(X,\mathcal{A})$. Since the scalar multiplication $S$ is continuous at $(1,f)$, there exists $\delta>0$ such that $S((1-\delta,1+\delta),U_{\mu,F}(f,\delta))\subseteqq U_{\mu,F}(f,1)$. Now $1+\frac{\delta}{2}\in (1-\delta,1+\delta)\implies S(1+\frac{\delta}{2},f)\in U_{\mu,F}(f,1)$ i.e., there exists $A\in\mathcal{A}$ with $\mu(A)<\infty$ such that $|(1+\frac{\delta}{2}).f-f|<1$ on $X\setminus A$. Thus $|f|<\frac{2}{\delta}$ on $X\setminus A\implies f\in L^\infty_F(X,\mathcal{A},\mu)$.
			\item Let $\mathcal{M}(X,\mathcal{A})=L^\infty_F(X,\mathcal{A},\mu)$, then from Lemma \ref{Lem6.8} it follows that $\mathcal{M}_{m_{\mu,F}}(X,\mathcal{A})=\mathcal{M}_{U_{\mu,F}}(X,\mathcal{A})$ and hence $\mathcal{M}_{m_{\mu,F}}(X,\mathcal{A})$ is a topological vector space. Proof of the converse part is similar to the above proof.  
		\end{enumerate}
	\end{proof}
	\begin{theorem}
		The following statements are equivalent:
		\begin{enumerate}
			\item $\mathcal{M}_{U^I_{\mu,F}}(X,\mathcal{A})$ is a topological vector space.
			\item $\mathcal{M}_{m^I_{\mu,F}}(X,\mathcal{A})$ is a topological vector space.
			\item $I=L^\infty_F(X,\mathcal{A},\mu)=\mathcal{M}(X,\mathcal{A})$.
		\end{enumerate}
	\end{theorem}
	\begin{proof}
		From Lemma \ref{Lem6.8} and Lemma \ref{Lem6.9}, it clearly follows that $(3)\implies(1)$ and $(3)\implies(2)$. Let $(1)$ hold i.e., $\mathcal{M}_{U^I_{\mu,F}}(X,\mathcal{A})$ is a topological vector space and therefore $L^\infty_F(X,\mathcal{A},\mu)=\mathcal{M}(X,\mathcal{A})$ (by Lemma \ref{Lem6.9}). If $f\in\mathcal{M}(X,\mathcal{A})$, then a similar approach as in Theorem \ref{Th2.5} implies $f\in I$. Therefore $I=L^\infty_F(X,\mathcal{A},\mu)=\mathcal{M}(X,\mathcal{A})$ i.e., $(1)\implies(3)$. Likewise $(2)\implies(3)$.
	\end{proof}
	\begin{definition}
		A set $A\in\mathcal{A}$ is called weakly $\mu$-bounded in $X$ if every $f\in\mathcal{M}(X,\mathcal{A})$ is weakly essentially bounded in $A$ i.e., for each $f\in \mathcal{M}(X,\mathcal{A})$, there exists a $\mu$-finite set $A_f\subseteqq A$ such that $f$ is bounded on $A\setminus A_f$. 
	\end{definition}
	It is clear that, every $\mu$-bounded subset of $X$ is weakly $\mu$-bounded. The next Theorem shows that when $U^I_{\mu,F}$-topology and $m^I_{\mu,F}$-topology are equal.
	\begin{theorem}
		For any ideal $I$ in $\mathcal{M}(X,\mathcal{A})$, if $X\setminus\bigcap Z[I]$ is a weakly $\mu$-bounded set in $X$, then $U^I_{\mu,F}$-topology $=m^I_{\mu,F}$-topology on $\mathcal{M}(X,\mathcal{A})$. When $I$ is finitely generated, then the converse is also true.
	\end{theorem}
	\begin{proof}
		The proof of the first part can be accomplished by closely following the proof of the first part in Theorem \ref{Th2.6}. For the converse part let $I=<f_i>_{i=1}^n$ be a finitely generated ideal in $\mathcal{M}(X,\mathcal{A})$ and $U^I_{\mu,F}$-topology $=m^I_{\mu,F}$-topology on $\mathcal{M}(X,\mathcal{A})$. It can be easily proved that $\bigcap Z[I]=\cap_{i=1}^nZ(f_i)\in\mathcal{A}$. Let $f\in \mathcal{M}(X,\mathcal{A})$. In $\mathcal{M}(X,\mathcal{A})$, consider the positive unit $u(x)=\begin{cases}
			\frac{1}{|f(x)|}&\text{ if }x\in X\setminus Z(f)\\
			1&\text{ if }x\in Z(f)
		\end{cases}$. Then $m_{\mu,F}(0,I,u)$ is an open neighbourhood of $0$ in the $m_{\mu,F}$-topology on $\mathcal{M}(X,\mathcal{A})$. Since $U^I_{\mu,F}$-topology $=m^I_{\mu,F}$-topology on $\mathcal{M}(X,\mathcal{A})$, there exists $\epsilon>0$ such that $U_{\mu,F}(0,I,\epsilon)\subseteqq m_{\mu,F}(0,I,u)$. Let $h(x)=\begin{cases}
		0&\text{ when }x\in \bigcap Z[I]\\
		\frac{\epsilon}{2}&\text{ otherwise}
	\end{cases}$. Then $Z(h)=\cap_{i=1}^n Z(f_i)\implies h\in I$. Also $0\leq h\leq\frac{\epsilon}{2}$. Thus $h\in U_{\mu,F}(0,I,\epsilon)\subseteqq m_{\mu,F}(0,I,u)\implies$ there exists $\mu$-finite set $A$ such that $h<u$ on $X\setminus A$. Now for all $x\in X\setminus A\cap X\setminus \bigcap Z[I]\cap X\setminus Z(f)$, $\frac{\epsilon}{2}<\frac{1}{|f(x)|}$ i.e., for all $x\in (X\setminus \bigcap Z[I])\setminus A$, $|f(x)|<\frac{2}{\epsilon}$. Let $A'=A\cap (X\setminus \bigcap Z[I])$. Then $A'\subset X\setminus \bigcap Z[I]$ and $\mu(A')<\infty$ (as $A$ is a $\mu$-finite set) and for all $x\in (X\setminus \bigcap Z[I])\setminus A'$, $|f(x)|<\frac{2}{\epsilon}$. Thus $f$ is weakly essentially bounded in $X\setminus \bigcap Z[I]$. This completes the proof.
	\end{proof}
	By choosing $I=\mathcal{M}(X,\mathcal{A})=<1>$, we have the following corollary:
	\begin{corollary}
		$U_{\mu,F}$-topology $=m_{\mu,F}$-topology on $\mathcal{M}(X,\mathcal{A})$ if and only if $\mathcal{M}(X,\mathcal{A})=L^\infty_F(X,\mathcal{A},\mu)$.
	\end{corollary}
	If we closely follow the arguments made in Section \ref{Sec3}, we can prove the following facts without any difficulty.
	\begin{theorem}\label{Th6.14}
		\hspace*{3cm}
		\begin{enumerate}
			\item Any additive subgroup of $\mathcal{M}(X,\mathcal{A})$ containing $I$ is clopen in $\mathcal{M}_{U^I_{\mu,F}}(X,\mathcal{A})$.
			\item $I\cap L^\infty_F(X,\mathcal{A},\mu)$ is a clopen subset of $\mathcal{M}_{U^I_{\mu,F}}(X,\mathcal{A})$.
			\item $I\cap L^\infty_F(X,\mathcal{A},\mu)$ is the component of $0$ in $\mathcal{M}_{U^I_{\mu,F}}(X,\mathcal{A})$.\label{Th6.14.3}
		\end{enumerate}
	\end{theorem}
	\begin{corollary}
		$L^\infty_F(X,\mathcal{A},\mu)$ is the component of $0$ in $\mathcal{M}_{U_{\mu,F}}(X,\mathcal{A})$.
	\end{corollary}
	\begin{definition}
		An ideal $I$ in $\mathcal{M}(X,\mathcal{A})$ is called weakly $\mu$-bounded if every $f\in I$ is weakly essentially bounded on $X$ i.e., $I\subseteqq L^\infty_F(X,\mathcal{A},\mu)$.
	\end{definition}
	\begin{notation}
		$L_{\psi,F}(X,\mathcal{A},\mu)=\{f\in\mathcal{M}(X,\mathcal{A}):~\forall g\in\mathcal{M}(X,\mathcal{A}), f.g\in L^\infty_F(X,\mathcal{A},\mu)\}$
	\end{notation}
	One can easily verify that $L_{\psi,F}(X,\mathcal{A},\mu)$ is the largest weakly $\mu$-bounded ideal in $\mathcal{M}(X,\mathcal{A})$. After some proper adjustment, the next two theorems follow from Theorem \ref{Th3.55} and Theorem \ref{Th3.5}.
	\begin{theorem}
		$L_{\psi,F}(X,\mathcal{A},\mu)=\{f\in\mathcal{M}(X,\mathcal{A}): X\setminus Z(f)\text{ is weakly }\mu\text{-bounded subset of }X\}$.
	\end{theorem}
	\begin{theorem}\label{Th6.18}
		The component of $0$ in $\mathcal{M}_{m^I_{\mu,F}}(X,\mathcal{A})$ is $I\cap L_{\psi,F}(X,\mathcal{A},\mu)$.
	\end{theorem}
	\begin{corollary}
		The component of $0$ in $\mathcal{M}_{m_{\mu,F}}(X,\mathcal{A})$ is $L_{\psi,F}(X,\mathcal{A},\mu)$, the largest weakly $\mu$-bounded ideal in $\mathcal{M}(X,\mathcal{A})$.
	\end{corollary}
	We recall that \begin{align*}
		\mathcal{M}_F(X,\mathcal{A})&=\{f\in\mathcal{M}(X,\mathcal{A}): \mu(X\setminus Z(f))<\infty\}\\
		\mathcal{M}_\infty(X,\mathcal{A})&=\{f\in\mathcal{M}(X,\mathcal{A}):~\forall \epsilon>0,\mu(\{x\in X:|f(x)|\geq\epsilon\})<\infty\}\\
		&=\{f\in\mathcal{M}(X,\mathcal{A}):~\forall n\in\mathbb{N},\mu(\{x\in X:|f(x)|\geq\frac{1}{n}\})<\infty\}\\
		\mathcal{M}_\psi(X,\mathcal{A})&=\{f\in\mathcal{M}(X,\mathcal{A}):~\forall g\in\mathcal{M}(X,\mathcal{A}),f.g\in\mathcal{M}_\infty(X,\mathcal{A})\}
	\end{align*}
	$\mathcal{M}_F(X,\mathcal{A})\subseteqq\mathcal{M}_\psi(X,\mathcal{A}) \subseteqq\mathcal{M}_\infty(X,\mathcal{A})$ and $\mathcal{M}_F(X,\mathcal{A}), \mathcal{M}_\psi(X,\mathcal{A})$ are ideals in $\mathcal{M}(X,\mathcal{A})$. If we closely observe Theorem \ref{Th4.0}, Theorem \ref{Th4.1} and Theorem \ref{Th4.4}, we get the following results.
	\begin{theorem}
		\hspace*{3cm}
		\begin{enumerate}
			\item $L^\infty_F$ is closed in the $U_{\mu,F}$-topology on $\mathcal{M}(X,\mathcal{A})$.
			\item $cl_{U_{\mu,F}}(\mathcal{M}_F(X,\mathcal{A}))=\mathcal{M}_\infty(X,\mathcal{A})$.
			\item $cl_{m_{\mu,F}}(\mathcal{M}_F(X,\mathcal{A}))=\mathcal{M}_\psi(X,\mathcal{A})$.
		\end{enumerate}
	\end{theorem}
	The notion of $Z_\mu$-ideals in $\mathcal{M}(X,\mathcal{A})$ has its analogue in the finite case also.
	\begin{definition}
		\hspace*{3cm}
		\begin{enumerate}
			\item Let $Y$ be a subset of $X$. Two functions $f,g\in \mathcal{M}(X,\mathcal{A})$ are called equal almost everywhere finitely $(\mu)$ on $Y$ if there exists a $\mu$-finite set $Z\subset Y$ such that $f\equiv g$ on $Y\setminus Z$ i.e., $\mu(\{y\in Y:f(y)\neq g(y)\})<\infty$. We write this case as $f\equiv_F g$ on $Y$.
			\item An ideal $I$ is said to be a $Z_{\mu,F}$-ideal if whenever $f\in\mathcal{M}(X,\mathcal{A})$ and $f\equiv_F g$ on $X$ for some $g\in I$, then $f\in I$. Equivalently, $I$ is a $Z_{\mu,F}$-ideal if for $f\in\mathcal{M}(X,\mathcal{A})$, $f\equiv_F 0$ on some $Z\in Z[I]$ implies $f\in I$.
		\end{enumerate}
	\end{definition}
	It is easy to check that, $\mathcal{L}^F(X,\mathcal{A},\mu), L_{\psi,F}(X,\mathcal{A},\mu)$ are $Z_{\mu,F}$-ideals in $\mathcal{L}^F(X,\mathcal{A},\mu)$. It is also immediate that, every $Z_{\mu,F}$-ideal in $\mathcal{M}(X,\mathcal{A})$ is also a $Z_\mu$-ideal. The next result closely follows from Theorem \ref{Th5.0}.
	\begin{theorem}
		An ideal $I$ in $\mathcal{M}(X,\mathcal{A})$ is a $Z_{\mu,F}$-ideal if and only if $I\supseteqq \mathcal{L}^F(X,\mathcal{A},\mu)$.
	\end{theorem}
	\begin{corollary}\label{Rem6.25}
		$\mathcal{L}^F(X,\mathcal{A},\mu)$ is the smallest $Z_{\mu,F}$-ideal in $\mathcal{M}(X,\mathcal{A})$.
	\end{corollary}
	The intersection of any collection of $Z_{\mu,F}$-ideals in $\mathcal{M}(X,\mathcal{A})$ is also a $Z_{\mu,F}$-ideal. Hence for any ideal $I$ in $\mathcal{M}(X,\mathcal{A})$, $I^e_F\equiv$ the intersection of all the $Z_{\mu,F}$-ideals in $\mathcal{M}(X,\mathcal{A})$ containing $I$, is the smallest $Z_{\mu,F}$-ideal containing $I$ in $\mathcal{M}(X,\mathcal{A})$. The next result gives us the formation of $I^e_F$, explicitly.
	\begin{theorem}
		Let $I$ be an ideal in $\mathcal{M}(X,\mathcal{A})$. If there exists an $f\in I$ such that $\mu(Z(f))<\infty$, then $I^e_F= \mathcal{M}(X,\mathcal{A})$. If for all $f\in I$, $\mu(Z(f))=\infty$, then $I^e_F=\{f\in\mathcal{M}(X,\mathcal{A}): \exists g\in I\text{ such that }\mu(Z(f)\triangle Z(g))<\infty\}$.
	\end{theorem}
	Proof of this Theorem can be made by closely following the proof of Theorem \ref{Th5.1} and so we omit the proof.\\
	Let the $\mathcal{A}$-filter $\mathscr{F'}=\{A\in\mathcal{A}:\mu(X\setminus A)<\infty\}$ be extended to an $\mathcal{A}$-ultrafilter $\mathscr{U'}$ and $M'$ be the maximal ideal in $\mathcal{M}(X,\mathcal{A})$ such that $Z[M']=\mathscr{U'}$. We can show as in Example \ref{Ex5.3} that $M'$ is a $Z_{\mu,F}$-ideal in $\mathcal{M}(X,\mathcal{A})$. Also the Example \ref{Ex5.4} shows that if $(X,\mathcal{A},\mu)$ is such a measure space that there exists $x\in X$ with $\{x\}\in\mathcal{A}$ and $\mu(\{x\})=0$, then $\mathcal{M}(X,\mathcal{A})$ contains at least a maximal ideal which is not a $Z_{\mu,F}$-ideal, given by $M_x=\{f\in \mathcal{M}(X,\mathcal{A}):f(x)=0\}$. Therefore a maximal ideal in $\mathcal{M}(X,\mathcal{A})$ may or may not be a $Z_{\mu,F}$-ideal. Let $\mathscr{C}$ be the collection of all $Z_{\mu,F}$-ideals in $\mathcal{M}(X,\mathcal{A})$. Then $\mathcal{L}^F(X,\mathcal{A},\mu)\in \mathscr{C}$. By using Zorn's lemma, we can find a maximal element in $\mathscr{C}$ containing $\mathcal{L}^F(X,\mathcal{A},\mu)$. We call the maximal element of $\mathscr{C}$ as the maximal $Z_{\mu,F}$-ideal in $\mathcal{M}(X,\mathcal{A})$. Later we show that every maximal $Z_{\mu,F}$-ideal is also a maximal ideal in $\mathcal{M}(X,\mathcal{A})$.
	\begin{theorem}
		\hspace*{3cm}
		\begin{enumerate}
			\item Any closed ideal in $\mathcal{M}_{m_{\mu,F}}(X,\mathcal{A})$ is a $Z_{\mu,F}$-ideal. 
			\item If closure of an ideal in the $U_{\mu,F}$-topology is an ideal in $\mathcal{M}(X,\mathcal{A})$, then that ideal is a $Z_{\mu,F}$-ideal.
		\end{enumerate}
	\end{theorem}
	\begin{corollary}
		For any ideal $I$ in $\mathcal{M}(X,\mathcal{A})$, $I\subseteqq I^e_F\subseteqq cl_{m_{\mu,F}}I\subseteqq cl_{U_{\mu,F}}I$.
	\end{corollary}
	\begin{theorem}
		For a proper ideal $I$ in $\mathcal{M}(X,\mathcal{A})$, the following statements are equivalent:
		\begin{enumerate}
			\item There exists $Z\in Z[I]$ such that $\mu(Z)<\infty$.
			\item $I$ is dense in $\mathcal{M}_{m_{\mu,F}}(X,\mathcal{A})$.
			\item $I$ is dense in $\mathcal{M}_{U_{\mu,F}}(X,\mathcal{A})$.
		\end{enumerate}
	\end{theorem}
	\begin{corollary}
		A proper $Z_{\mu,F}$-ideal in $\mathcal{M}(X,\mathcal{A})$ is never dense in $\mathcal{M}_{m_{\mu,F}}(X,\mathcal{A})$ or in $\mathcal{M}_{U_{\mu,F}}(X,\mathcal{A})$.
	\end{corollary}
	\begin{notation}
		For any ideal $I$ in $\mathcal{M}(X,\mathcal{A})$, let $I_{\mu,F}=\{f\in\mathcal{M}(X,\mathcal{A}):~\forall g\in\mathcal{M}(X,\mathcal{A}), f.g\in cl_{U_{\mu,F}}I\}$.
	\end{notation}
	\begin{theorem}
		Let $I$ be an ideal in $\mathcal{M}(X,\mathcal{A})$. Then,
		\begin{enumerate}
			\item $I_{\mu,F}$ is an ideal in $\mathcal{M}(X,\mathcal{A})$ containing $I$.
			\item $I_{\mu,F}$ is the largest ideal containing $I$ and contained in $cl_{U_{\mu,F}}I$.
			\item $I_{\mu,F}$ is a closed subset of $\mathcal{M}(X,\mathcal{A})$ in the $m_{\mu,F}$-topology.
			\item $I_{\mu,F}=cl_{m_{\mu,F}}I$.
			\item $I_{\mu,F}$ is a $Z_{\mu,F}$-ideal.
			\item $I^e_F=I_{\mu,F}$.
			\item $I$ is a $Z_{\mu,F}$-ideal if and only if $I=I_{\mu,F}$.
		\end{enumerate}
	\end{theorem}
	\begin{corollary}\label{Cor6.29}
		An ideal in $\mathcal{M}(X,\mathcal{A})$ is closed in the $m_{\mu,F}$-topology if and only if it is a $Z_{\mu,F}$-ideal.
	\end{corollary}
	\begin{theorem}
		$\{I\subset\mathcal{M}(X,\mathcal{A}):\mathcal{L}^F(X,\mathcal{A},\mu) \subseteqq I\text{ and }I\text{ is an ideal in }\mathcal{M}(X,\mathcal{A})\}$ is precisely the set of all ideals in $\mathcal{M}(X,\mathcal{A})$ which are closed in the $m_{\mu,F}$-topology.
	\end{theorem}
	\begin{corollary}
		For an ideal $I$ in $\mathcal{M}(X,\mathcal{A})$, $cl_{U_{\mu,F}}I$ is an ideal in $\mathcal{M}(X,\mathcal{A})$ if and only if $cl_{U_{\mu,F}}I=cl_{m_{\mu,F}}I$.
	\end{corollary}
	\begin{theorem}\label{Th6.31}
		Let $M$ be a maximal ideal in $\mathcal{M}(X,\mathcal{A})$. 
		\begin{enumerate}
			\item $M$ is a $Z_{\mu,F}$-ideal if and only if $\mu(Z(f))=\infty$ for all $f\in M$.
			\item $M$ is closed in the $m_{\mu,F}$-topology if and only if $\mu(Z(f))=\infty$ for all $f\in M$.\label{Th6.31.3}
		\end{enumerate}
	\end{theorem}
	\begin{corollary}
		A maximal $Z_{\mu,F}$-ideal in $\mathcal{M}(X,\mathcal{A})$ is also a maximal ideal.
	\end{corollary}
	\begin{theorem}
		Let $(X,\mathcal{A},\mu)$ be a measure space where $\{x\}\in\mathcal{A}$ for all $x\in X$. Then the following two statements are equivalent:
		\begin{enumerate}
			\item Each $Z_{\mu,F}$-ideal in $\mathcal{M}(X,\mathcal{A})$ is free.
			\item $\mu(\{x\})<\infty$ for all $x\in X$. 
		\end{enumerate}
	\end{theorem}
	All these facts can be established from Section \ref{Sec5} after some necessary modification. Therefore we omit the proofs of these results.
	\section{Illustrative Examples}\label{Sec7}
	\counterwithin{theorem}{subsection}
	In this section, we record a few basic facts concerning the $U_\mu^I$-topology, $m_\mu^I$-topology, $U_{\mu,F}^I$-topology and $m_{\mu,F}^I$-topology on $\mathcal{M}(X,\mathcal{A})$ for some well-known measure spaces $(X,\mathcal{A},\mu)$.
	\subsection{Counting Measure}
	Let $X$ be a non-empty set, $\mathcal{A}=\mathscr{P}(X)$, the power set of $X$ and $\mu$ be the counting measure on $X$ defined as follows: $\mu(E)=\begin{cases}
		|E|\text{ if }E\text{ is finite}\\
		\infty\text{ otherwise}
	\end{cases}$\\ The following results can be proved with any difficulties.
	\begin{enumerate}
		\item $\mathcal{M}(X,\mathcal{A})=\mathbb{R}^X$ i.e., every real-valued function on $X$ is measurable.
		\item An $f\in\mathcal{M}(X,\mathcal{A})$ is essentially bounded on $X$ if and only if $f$ is bounded on $X$.
		\item $L^\infty(X,\mathcal{A},\mu) =\mathcal{M}(X,\mathcal{A})$ if and only if $X$ is finite. 
		\item $\mathcal{M}_{U^I_\mu}(X,\mathcal{A})$ (and $\mathcal{M}_{m^I_\mu}(X,\mathcal{A})$) is a topological vector space if and only if $X$ is finite and $I=\mathcal{M}(X,\mathcal{A})$.
		\item A subset $Y$ of $X$ is $\mu$-bounded if and only if $Y$ is a finite set.
		\item If $X\setminus \bigcap Z[I]$ is a finite set, then $\mathcal{M}_{U^I_\mu}(X,\mathcal{A})= \mathcal{M}_{m^I_\mu}(X,\mathcal{A})$. Conversely if $I$ is countably generated and $\mathcal{M}_{U^I_\mu}(X,\mathcal{A})= \mathcal{M}_{m^I_\mu}(X,\mathcal{A})$, then $X\setminus \bigcap Z[I]$ is finite.
		\item The component of $0$ in the $U^I_\mu$-topology is $\{f\in I:f\text{ is bounded on }X\}$.
		\item $L_\psi(X,\mathcal{A},\mu)=\{f\in\mathcal{M}(X,\mathcal{A}):X\setminus Z(f)\text{ is finite}\}$.
		\item The component of $0$ in the $m^I_\mu$-topology is $\{f\in I:X\setminus Z(f)\text{ is finite}\}$.
		\item For $f,g\in\mathcal{M}(X,\mathcal{A})$, $f\equiv g$ almost everywhere $(\mu)$ on $X$ if and only if $f\equiv g$ on $X$.
		\item Every ideal in $\mathcal{M}(X,\mathcal{A})$ is a $Z_\mu$-ideal.
	\end{enumerate}
	We now consider $X$ to be an infinite set. Then the following arguments hold.
	\begin{enumerate}
		\item An $f\in\mathcal{M}(X,\mathcal{A})$ is weakly essentially bounded if and only if $f$ is bounded except on a finite subset of $X$ if and only if $f$ is bounded on $X$.
		\item $L^\infty_F(X,\mathcal{A},\mu)=\{f\in\mathcal{M}(X,\mathcal{A}):f\text{ is bounded on }X\}=L^\infty(X,\mathcal{A},\mu)$.
		\item $\mathcal{M}_{U^I_{\mu,F}}(X,\mathcal{A})$ (and $\mathcal{M}_{m^I_{\mu,F}}(X,\mathcal{A})$) can never be a topological vector space.
		\item Weakly $\mu$-bounded sets and $\mu$-bounded sets are equal i.e., a subset $A$ of $X$ is weakly $\mu$-bounded if and only if $A$ is finite.
		\item If $X\setminus\bigcap Z[I]$ is a finite set for some ideal $I$ in $\mathcal{M}(X,\mathcal{A})$, then $U^I_{\mu,F}$-topology $=m^I_{\mu,F}$-topology on $\mathcal{M}(X,\mathcal{A})$. Conversely if $I$ is finitely generated ideal in $\mathcal{M}(X,\mathcal{A})$ and $U^I_{\mu,F}$-topology $=m^I_{\mu,F}$-topology on $\mathcal{M}(X,\mathcal{A})$, then $X\setminus\bigcap Z[I]$ is a finite subset of $X$.
		\item The component of $0$ in the $U_{\mu,F}^I$-topology on $\mathcal{M}(X,\mathcal{A})$ is $\{f\in I:f\text{ is bounded on }X\}$.
		\item $L_{\psi,F}(X,\mathcal{A},\mu)=\{f\in\mathcal{M}(X,\mathcal{A}):X\setminus Z(f)\text{ is a finite subset of }X\}$.
		\item The component of $0$ in the $m_{\mu,F}^I$-topology on $\mathcal{M}(X,\mathcal{A})$ is $\{f\in I:X\setminus Z(f)\text{ is a finite subset of }X\}$.
		\item Let $f,g\in\mathcal{M}(X,\mathcal{A})$. Then $f\equiv_Fg$ on $X$ if and only if $f\equiv g$ except on a finite subset of $X$.
		\item An ideal $I$ in $\mathcal{M}(X,\mathcal{A})$ is dense in both $U_{\mu,F}$ and $m_{\mu,F}$ topologies on $\mathcal{M}(X,\mathcal{A})$ if and only if $I$ contains an $f$ such that $Z(f)$ is a finite subset of $X$.
		\item The smallest $Z_{\mu,F}$-ideal in $\mathcal{M}(X,\mathcal{A})$ is given be $\mathcal{L}^F(X,\mathcal{A},\mu)=\{f\in\mathcal{M}(X,\mathcal{A}) :X\setminus Z(f)\text{ is a finite subset of }X\}$.
	\end{enumerate}
	\subsection{Dirac Measure}
	Let $X$ be a non-empty set, $\mathcal{A}=\mathscr{P}(X)$, the power set of $X$ and $x_0\in X$. The Dirac measure at $x_0$ is defined by $\delta_{x_0}(A)=\begin{cases}
		0&\text{ if }x_0\notin A\\
		1&\text{ if }x_0\in A
	\end{cases}$\\ Then the following statements are easy to show.
	\begin{enumerate}
		\item $\mathcal{M}(X,\mathcal{A})=\mathbb{R}^X$.
		\item Every $f\in\mathcal{M}(X,\mathcal{A})$ is essentially bounded.
		\item $L^\infty(X,\mathcal{A},\delta_{x_0})=\mathcal{M}(X,\mathcal{A})$.
		\item $\mathcal{M}_{U^I_{\delta_{x_0}}}(X,\mathcal{A}), \mathcal{M}_{m^I_{\delta_{x_0}}}(X,\mathcal{A})$ are topological vector spaces if and only if $I=\mathcal{M}(X,\mathcal{A})$.
		\item Every subset $Y$ of $X$ is $\delta_{x_0}$-bounded.
		\item $\mathcal{M}_{U^I_{\delta_{x_0}}}(X,\mathcal{A})= \mathcal{M}_{m^I_{\delta_{x_0}}}(X,\mathcal{A})$ for every ideal $I$ of $\mathcal{M}(X,\mathcal{A})$.
		\item For every ideal $I$ in $\mathcal{M}(X,\mathcal{A})$, the component of $0$ in the $U^I_{\delta_{x_0}}$-topology is $I$.
		\item $L_\psi(X,\mathcal{A},\delta_{x_0})=\mathcal{M}(X,\mathcal{A})$.
		\item For every ideal $I$ in $\mathcal{M}(X,\mathcal{A})$, the component of $0$ in the $m^I_{\delta_{x_0}}$-topology is $I$.
		\item Two functions $f,g\in\mathcal{M}(X,\mathcal{A})$ are equal almost everywhere $(\delta_{x_0})$ on $X$ if and only if $f(x_0)=g(x_0)$.
		\item An $f\in\mathcal{M}(X,\mathcal{A})$ is a $\delta_{x_0}$-unit if and only if $f(x_0)\neq 0$.
		\item $\mathcal{L}^0(X,\mathcal{A},\delta_{x_0})= \{f\in\mathcal{M}(X,\mathcal{A}):f(x_0)=0\}$.
		\item An ideal $I$ in $\mathcal{M}(X,\mathcal{A})$ is dense in both $U_{\delta_{x_0}}$ and $m_{\delta_{x_0}}$ topologies if and only if $x_0\notin\bigcap Z[I]$.
	\end{enumerate}
	\begin{theorem}
		$M_{x_0}=\{f\in\mathcal{M}(X,\mathcal{A}):f(x_0)=0\}$ is the only proper $Z_{\delta_{x_0}}$-ideal in $\mathcal{M}(X,\mathcal{A})$.
	\end{theorem}
	\begin{proof}
		We already have, $\mathcal{L}^0(X,\mathcal{A},\delta_{x_0})=M_{x_0}$. Therefore by Corollary \ref{Rem5.2}, $M_{x_0}$ is the smallest $Z_{\delta_{x_0}}$-ideal. Since $M_{x_0}$ is a maximal ideal, $M_{x_0}$ is the only proper $Z_{\delta_{x_0}}$-ideal in $\mathcal{M}(X,\mathcal{A})$.
	\end{proof}
	\begin{corollary}
		For any ideal $I$ in $\mathcal{M}(X,\mathcal{A})$, either $I\subset M_{x_0}$ and in this case $cl_{m_{\delta_{x_0}}}I=M_{x_0}$ or $I$ is dense in both $U_{\delta_{x_0}}$ and $m_{\delta_{x_0}}$ topologies on $\mathcal{M}(X,\mathcal{A})$.
	\end{corollary}
	\subsection{Co-countable Measure}\label{CoCu}
	Let $X$ be an uncountable set, $\mathcal{A}=\{A\subseteqq X:\text{ either }A\text{ or }X\setminus A\text{ is countable}\}$. The co-countable measure is defined by $\rho(A)=\begin{cases}
		0&\text{ if }A\text{ is countable}\\
		\infty&\text{ if }X\setminus A\text{ is countable}
	\end{cases}$.
	\begin{theorem}
		$\mathcal{M}(X,\mathcal{A})=\{f\in\mathbb{R}^X:f\text{ is constant except on a countable subset of }X\}$.
	\end{theorem}
	\begin{proof}
		If $f:X\to\mathbb{R}$ is a function such that it is constant except on a countable subset of $X$, then it is indeed a measurable function on $X$. Now let $f\in\mathcal{M}(X,\mathcal{A})$. For each $r\in\mathbb{R}$, let $A_r=\{x\in X:f(x)=r\}=f^{-1}(\{r\})\in\mathcal{A}$. If $X\setminus A_r$ is countable for some $r\in\mathbb{R}$, then there is nothing to prove. If possible let $X\setminus A_r$ be uncountable for all $r\in\mathbb{R}$. As $A_r\in\mathcal{A}$ for all $r\in\mathbb{R}$, this implies that $A_r$ is countable for each $r\in\mathbb{R}$. Let $S=\{r\in\mathbb{R}: A_r\neq\emptyset\}$. Since $X$ is uncountable, $X=\bigsqcup\limits_{r\in \mathbb{R}}A_r$ and $A_r$ is countable, it follows that $S$ is uncountable and so we can find an $s\in S$ such that $S_1=\{r\in S:r<s\}$, $S_2=\{r\in S:r\geq s\}$ are uncountable subsets of $S$. Now $f^{-1}(-\infty, s)=\bigsqcup\limits_{r\in S_1}A_r\implies f^{-1}(-\infty, s)$ is an uncountable set and $X\setminus f^{-1}(-\infty, s)=f^{-1}[s,\infty)=\bigsqcup\limits_{r\in S_2}A_r\implies X\setminus f^{-1}(-\infty, s)$ is an uncountable set. Again since $f\in \mathcal{M}(X,\mathcal{A})$, either $f^{-1}(-\infty, s)$ or $X\setminus f^{-1}(-\infty, s)$ is countable -- a contradiction.
	\end{proof}
	The following results are easy to establish.
	\begin{enumerate}
		\item Every $f\in\mathcal{M}(X,\mathcal{A})$ is essentially bounded.
		\item $L^\infty(X,\mathcal{A},\rho)=\mathcal{M}(X,\mathcal{A})$.
		\item $\mathcal{M}_{U^I_\rho}(X,\mathcal{A})$ and $\mathcal{M}_{m^I_\rho}(X,\mathcal{A})$ are always topological vector spaces if and only if $I=\mathcal{M}(X,\mathcal{A})$.
		\item Every set $A\in\mathcal{A}$ is $\rho$-bounded.
		\item For any ideal $I$ in $\mathcal{M}(X,\mathcal{A})$, $U^I_\rho$-topology $=m^I_\rho$-topology on $\mathcal{M}(X,\mathcal{A})$.
		\item For every ideal $I$ in $\mathcal{M}(X,\mathcal{A})$, the component of $0$ in the $U^I_\rho$-topology is $I$.
		\item $L_\psi(X,\mathcal{A},\rho)=\mathcal{M}(X,\mathcal{A})$.
		\item For every ideal $I$ in $\mathcal{M}(X,\mathcal{A})$, the component of $0$ in the $m^I_\rho$-topology is $I$.
		\item For $f,g\in\mathcal{M}(X,\mathcal{A})$, $f\equiv g$ almost everywhere $(\rho)$ on $X$ if and only if $f,g$ are equal except on a countable subset of $X$.
		\item $f\in\mathcal{M}(X,\mathcal{A})$ is a $\rho$-unit if and only if $f$ vanishes atmost countably many points of $X$.
		\item An ideal $I$ in $\mathcal{M}(X,\mathcal{A})$ is dense in both $U_\rho$ and $m_\rho$ topologies if and only if $Z(f)$ is countable for some $f\in I$.
		\item $\mathcal{L}^0(X,\mathcal{A},\rho)= \{g\in\mathcal{M}(X,\mathcal{A}):X\setminus Z(g)\text{ is countable}\}=\mathcal{L}^F(X,\mathcal{A},\rho)$.
	\end{enumerate}
	\begin{theorem}\label{Th7.4}
		$\mathcal{L}^0(X,\mathcal{A},\rho)$ is the only proper $Z_\rho$-ideal in $\mathcal{M}(X,\mathcal{A})$.
	\end{theorem}
	\begin{proof}
		From Corollary \ref{Rem5.2}, it follows that $\mathcal{L}^0(X,\mathcal{A},\rho)$ is the smallest $Z_\rho$-ideal in $\mathcal{M}(X,\mathcal{A})$. Let $I$ be a proper $Z_\rho$-ideal in $\mathcal{M}(X,\mathcal{A})$. Earlier we observed that if $I$ contains an $f$ such that $Z(f)$ is countable, then $I$ is dense in the $m_\rho$-topology on $\mathcal{M}(X,\mathcal{A})$. Since $I$ is a proper $Z_\rho$-ideal, $Z(f)$ is uncountable for all $f\in I\implies X\setminus Z(f)$ is countable for all $f\in I$. Therefore, $I\subseteqq \mathcal{L}^0(X,\mathcal{A},\rho)$.  Consequently $I= \mathcal{L}^0(X,\mathcal{A},\rho)$.
	\end{proof}
	\begin{corollary}
		For any ideal $I$ in $\mathcal{M}(X,\mathcal{A})$, either $I\subset \mathcal{L}^0(X,\mathcal{A},\rho)$ and in this case $cl_{m_\rho}I=\mathcal{L}^0(X,\mathcal{A},\rho)$ or $I$ is dense in both $U_\rho$ and $m_\rho$ topologies on $\mathcal{M}(X,\mathcal{A})$.
	\end{corollary}
	\begin{theorem}\label{Th7.5}
		For every ideal $I$ in $\mathcal{M}(X,\mathcal{A})$, $U_{\rho,F}^I$-topology $=U_\rho^I$-topology and $m_{\rho,F}^I$-topology $=m_\rho^I$-topology on $\mathcal{M}(X,\mathcal{A})$.
	\end{theorem}
	\begin{proof}
		Since for $A\in\mathcal{A}$, $\rho(A)<\infty\implies\rho(A)=0$, then $U_\rho(f,I,\epsilon)=U_{\rho,F}(f,I,\epsilon)$ for all $f\in\mathcal{M}(X,\mathcal{A})$ and $\epsilon>0$. Therefore $U_{\rho,F}^I$-topology $=U_\rho^I$-topology on $\mathcal{M}(X,\mathcal{A})$ and similarly $m_{\rho,F}^I$-topology $=m_\rho^I$-topology on $\mathcal{M}(X,\mathcal{A})$.
	\end{proof}
	Since these pair of topologies are same on $\mathcal{M}(X,\mathcal{A})$, all the results hold for the $U_\rho^I$-topology and $m_\rho^I$-topology are also true for $U_{\rho,F}^I$-topology and $m_{\rho,F}^I$-topology on $\mathcal{M}(X,\mathcal{A})$. Likewise we have the following result.
	\begin{corollary}
		$\mathcal{L}^0(X,\mathcal{A},\rho)$ is the only proper $Z_{\rho,F}$-ideal in $\mathcal{M}(X,\mathcal{A})$.
	\end{corollary}
	\subsection{Lebesgue Measure on $\mathbb{R}$}
	Let $X=\mathbb{R}$, $\mathcal{A}$ be the set of all Lebesgue measurable sets on $\mathbb{R}$ and $\mu$ be the Lebesgue measure $X$.
	\begin{theorem}\label{Th7.7}
		Let $A\in\mathcal{A}$. Then $A$ is $\mu$-bounded if and only if $\mu(A)=0$.
	\end{theorem}
	\begin{proof}
		It is trivial that if $\mu(A)=0$, then $A$ is $\mu$-bounded. Conversely let $\mu(A)>0$. Since it is a non-atomic measure, we can find a sequence of pairwise disjoint measurable subsets $\{A_n\}_{n=1}^\infty$ such that $\mu(A_n)>0$ for all $n\in\mathbb{N}$ and $A=\bigsqcup\limits_{n=1}^\infty A_n$. Let us define $f:X\to\mathbb{R}$ as follows: $f(x)=\begin{cases}
			n&\text{ if }x\in A_n\\
			0&\text{ Otherwise}
		\end{cases}$. Then for each $n\in\mathbb{N}$, $\{x\in A:f(x)\geq n\}=\bigsqcup\limits_{p=n}^\infty A_p\implies \mu(\{x\in A:f(x)\geq n\})>0$. Consequently, $f$ is not essentially bounded on $A$.
	\end{proof}
	We can generalize this theorem for any non-atomic measure space.
	\begin{theorem}\label{Th7.3}
		Let $\mu$ be a non-atomic measure on a measurable space $(X,\mathcal{A})$. Then $A\in\mathcal{A}$ is a $\mu$-bounded subset of $X$ if and only if $\mu(A)=0$.
	\end{theorem}
	We now show that a $\mu$-bounded subset is of zero measure characterizes the non-atomicity of the measure.
	\begin{lemma}\label{Lem7.4}
		Any atom $A$ in a measure space $(X,\mathcal{A},\mu)$ is a $\mu$-bounded subset of $X$.
	\end{lemma}
	\begin{proof}
		Let $f\in\mathcal{M}(X,\mathcal{A})$ and $A_n(f)=\{x\in X:|f(x)|\geq n\}$, for each $n\in\mathbb{N}$. Since $A$ is an atom, either $\mu(A\cap A_n(f))=0$ or $\mu(A\cap X\setminus A_n(f))=0$ for each $n\in\mathbb{N}$. Suppose that $\mu(A\cap A_n(f))=0$ for some $n\in\mathbb{N}$. Then $|f|<n$ on $A\setminus (A\cap A_n(f))$ i.e., $f$ is essentially bounded on $A$. If possible let $\mu(A\cap A_n(f))>0$ for all $n\in\mathbb{N}$. Then $\mu(A\cap X\setminus A_n(f))=0$ for all $n\in\mathbb{N}$. Now if $x\in A$, then by the Archimedean property of $\mathbb{R}$ there exists an $n\in\mathbb{N}$ such that $|f(x)|<n$ i.e., $x\in X\setminus A_n(f)$. Thus $A=\bigcup\limits_{n\in\mathbb{N}}(A\cap X\setminus A_n(f))\implies\mu(A)=0$ -- contradicts that $A$ is an atom.
	\end{proof}
	\begin{theorem}\label{Th7.10}
		In a measurable space $(X,\mathcal{A})$, $\mu$ is a non-atomic measure if and only if every $\mu$-bounded subset is of measure zero.
	\end{theorem}
	\begin{proof}
		It follows from Theorem \ref{Th7.3} that if $\mu$ is a non-atomic measure, then a set $A\in\mathcal{A}$ is $\mu$-bounded if and only if $\mu(A)=0$. Conversely if $\mu$ is not a non-atomic measure, then $X$ contains an atom, say $A$. Then by Lemma \ref{Lem7.4}, $A$ is $\mu$-bounded with $\mu(A)>0$.
	\end{proof}
	\begin{corollary}
		A measure space $(X,\mathcal{A},\mu)$ is non-atomic if and only if $L_\psi(X,\mathcal{A},\mu)=\{f\in\mathcal{M}(X,\mathcal{A}):\mu(X\setminus Z(f))=0\}$ i.e., $L_\psi(X,\mathcal{A},\mu)= \mathcal{L}^0(X,\mathcal{A},\mu)$.
	\end{corollary}
	Let $(X,\mathcal{A},\mu)$ be a non-atomic measure space. Then the following statements are true.
	\begin{enumerate}
		\item $L^\infty(X,\mathcal{A},\mu)\neq\mathcal{M}(X,\mathcal{A})$.
		\item $\mathcal{M}_{U^I_\mu}(X,\mathcal{A})$ (and $\mathcal{M}_{m^I_\mu}(X,\mathcal{A})$) can never be a topological vector space.
		\item The component of $0$ in the $m^I_\mu$-topology is $\{f\in I:\mu(X\setminus Z(f))=0\}$.
		\item If $\mu(X\setminus\bigcap Z[I])=0$ for some ideal $I$ in $\mathcal{M}(X,\mathcal{A})$, then $U_\mu^I$-topology $=m_\mu^I$-topology on $\mathcal{M}(X,\mathcal{A})$.
	\end{enumerate}
	In particular, all these results also hold for the Lebesgue measure space $(\mathbb{R},\mathcal{A},\mu)$ on $\mathbb{R}$. We now find out which are the weakly $\mu$-bounded subsets of $\mathbb{R}$.
	\begin{lemma}\label{Lem7.11}
		Let $(\mathbb{R},\mathcal{A},\mu)$ be the Lebesgue measure space and $A\in\mathcal{A}$ be such that $\mu(A)=\infty$.
		\begin{enumerate}
			\item For every $n\in\mathbb{N}$, there exists $A_n\subset A$ such that $\mu(A_n)=n$.
			\item There exists a sequence $\{A_n\}_{n=1}^\infty$ of pairwise disjoints measurable subsets of $A$ such that $\mu(A_n)=n$ for all $n\in\mathbb{N}$.
		\end{enumerate}
	\end{lemma}
	\begin{proof}
		Since $\mu(A)=\infty$, then either $\mu(A\cap[0,\infty])=\infty$ or $\mu(A\cap[-\infty,0])=\infty$. Without loss of generality let $\mu(A\cap[0,\infty])=\infty$.
		\begin{enumerate}
			\item Fix an $n\in\mathbb{N}$. Let us define $f:\mathbb{R}\to[0,\infty]$ as follows: $$f(x)=\begin{cases}
				0&\text{ if }x\leq 0\\
				\mu(A\cap[0,\infty])&\text{ if }x\geq 0
			\end{cases}$$
			Then $f$ is continuous on $\mathbb{R}$ because for all $x,y\geq 0$, $|f(x)-f(y)|<|x-y|$. Also $\lim\limits_{x\to 0}f(x)=0$ and $\lim\limits_{x\to\infty}=\mu(A\cap[0,\infty])=\infty$. Thus there exists an $x_n\in \mathbb{R}$ such that $f(x_n)=n$. Choose $A_n=A\cap[0,x_n]$. Then $\mu(A_n)=f(x_n)=n$.
			\item By the above result, there exists $A_1\subset A$ such that $\mu(A_1)=1$. Now $\mu(A\setminus A_1)=\infty$. Again using the above result, we can find $A_2\subset A\setminus A_1$ such that $\mu(A_2)$. Clearly $A_1\cap A_2=\emptyset$. Continuing this process, we get a sequence $\{A_n\}_{n=1}^\infty$ of pairwise disjoints measurable subsets of $A$ such that $\mu(A_n)=n$ for all $n\in\mathbb{N}$.
		\end{enumerate}
	\end{proof}
	\begin{theorem}\label{Th7.4.6}
		Let $(\mathbb{R},\mathcal{A},\mu)$ be the Lebesgue measure space and $A\in\mathcal{A}$. Then $A$ is a weakly $\mu$-bounded subset of $\mathbb{R}$ if and only if $\mu(A)<\infty$.
	\end{theorem}
	\begin{proof}
		If $\mu(A)<\infty$, then $A$ is trivially weakly $\mu$-bounded. Conversely let $\mu(A)=\infty$. Then by Lemma \ref{Lem7.11}, we get a sequence $\{A_n\}_{n=1}^\infty$ of pairwise disjoints measurable subsets of $A$ such that $\mu(A_n)=n$ for all $n\in\mathbb{N}$. Let us define $f:\mathbb{R}\to \mathbb{R}$ as follows: $f(x)=\begin{cases}
			n&\text{ if }x\in A_n\\
			0&\text{ Otherwise}
		\end{cases}$. Then for each $n\in\mathbb{N}$, $\{x\in A:f(x)\geq n\}=\bigsqcup\limits_{p=n}^\infty A_p\implies \mu(\{x\in A:f(x)\geq n\})=\infty$. Consequently, $f$ is not weakly essentially bounded on $\bigsqcup\limits_{n=1}^\infty A_n\subseteqq A$. 
	\end{proof}
	We state the following theorem due to Sierpi\'nski \cite{Sierpinski}.
	\begin{theorem}\label{Th7.4.7}
		In a non-atomic measure space $(X,\mathcal{A},\mu)$, for every positive measurable set $A\in\mathcal{A}$, we can find $B_\alpha(\subseteqq A)\in\mathcal{A}$ such that $\mu(B_\alpha)=\alpha$ for each $\alpha\in[0,\mu(A)]$.
	\end{theorem}
	In view of this Theorem \ref{Th7.4.7} and Lemma \ref{Lem7.11}, we can extend Theorem \ref{Th7.4.6} to any non-atomic measure space.
	\begin{theorem}\label{Th7.4.8}
		In a non-atomic measure space $(X,\mathcal{A},\mu)$, a measurable set $A$ is weakly $\mu$-bounded if and only if $\mu(A)<\infty$.
	\end{theorem}
	Furthermore, the weak $\mu$-boundedness of a measurable set can determine the atomicity of the measure.
	\begin{theorem}
		A measure space $(X,\mathcal{A},\mu)$ is non-atomic if and only if every weakly $\mu$-bounded set is $\mu$-finite.
	\end{theorem}
	\begin{proof}
		It follows from Theorem \ref{Th7.4.8} that in a non-atomic measure space $(X,\mathcal{A},\mu)$, a measurable set $A$ is $\mu$-bounded if and only if $\mu(A)<\infty$. Conversely suppose that $(X,\mathcal{A},\mu)$ contain an atom, say $A$. Then by Lemma \ref{Lem7.4}, $A$ is $\mu$-bounded subset of $X$. Consequently, $A$ is weakly $\mu$-bounded with $\mu(A)>0$.
	\end{proof}
	\begin{corollary}
		A measure space $(X,\mathcal{A},\mu)$ is non-atomic if and only if $L_{\psi, F}(X,\mathcal{A},\mu)=\{f\in\mathcal{M}(X,\mathcal{A}) :\mu(X\setminus Z(f))<\infty\}$ i.e., $L_{\psi, F}(X,\mathcal{A},\mu)= \mathcal{L}^F(X,\mathcal{A},\mu)$.
	\end{corollary}
	For a non-atomic measure space $(X,\mathcal{A},\mu)$, we have the following statements:
	\begin{enumerate}
		\item $L^\infty_F(X,\mathcal{A},\mu)\neq\mathcal{M}(X,\mathcal{A})$.
		\item $\mathcal{M}(X,\mathcal{A})$ equipped with $U_{\mu,F}^I$-topology and $m_{\mu,F}^I$-topology are not topological vector spaces.
		\item If $\mu(X\setminus\bigcap Z[I])<\infty$ for some ideal $I$ in $\mathcal{M}(X,\mathcal{A})$, then $U_{\mu,F}^I$-topology $=m_{\mu,F}^I$-topology on $\mathcal{M}(X,\mathcal{A})$.
		\item The component of $0$ in the $m_{\mu,F}^I$-topology on $\mathcal{M}(X,\mathcal{A})$ is $\{f\in I:\mu(X\setminus Z(f))<\infty\}$.
	\end{enumerate}
		

\begin{thebibliography}{1}
		\bibitem{Acharyya1}
		Acharyya, S. K., Bag, S. and Sack, J., \emph{Ideals of rings and intermediate rings of measurable functions}, J. of Alg. and Its App., Vol. 19(2) (2020), \href{https://doi.org/10.1142/S0219498820500383}{2050038}.
				
		\bibitem{Acharyya2}
		Acharyya, S., Acharyya, S. K., Bag, S. and Sack, J., \emph{Recent progress in rings and subrings of real valued measurable functions}, Quaes. Math., Vol. 43(7) (2020), \href{https://doi.org/10.2989/16073606.2019.1585395}{959--973}.	
			
		\bibitem{Acharyya3}
		Acharyya, S., Acharyya, S. K., Bharati, R. and Deb Ray, A., \emph{Some algebraic and topological properties of rings of measurable functions}, Hous. J. of Math., Vol 47(3) (2021), 633--657.
		
		\bibitem{Acharyya4}
		Acharyya, S., Bharati, R., Deb Ray, A. and Acharyya, S. K., \emph{A generalization of $m$-topology and $U$-topology on rings of measurable functions}, \textbf{communicated}.
				
		\bibitem{Azarpanah5}
		Azarpanah, F., Manshoor, F. and Mohamadian, R., \emph{A Generalization of the $m$-Topology on $C(X)$ Finer than the $m$-Topology}, Filomat, Vol. 31(8) (2017), \href{https://doi.org/10.2298/FIL1708509A}{2509--2515}.
		
		\bibitem{Bruckner6}
		Bruckner, A. M., Bruckner, J. B. and Thomson, B. S., \emph{Real Analysis}, ClassicalRealAnalysis.com, 2008, xiv 656 pp., \href{http://classicalrealanalysis.info/com/documents/BBT-AlllChapters-Landscape.pdf}{ISBN 1434844129}.
		
		\bibitem{Gillman7}
		Gillman, L. and Jerison, M., \emph{Rings of Continuous Functions}, D. Van Nos. Com., New York, \href{https://doi.org/10.1007/978-1-4615-7819-2}{1960}.
		
		\bibitem{Johnson}
		Johnson, R. A., \emph{Atomic and nonatomic measures}, Proc. Amer. Math. Soc., Vol 25 (1970), \href{https://doi.org/10.1090/S0002-9939-1970-0279266-8 
		}{650--655}. 
		
		\bibitem{Mandelker}
		Mandelker, M., \emph{Supports of continuous functions}, Trans. Amer. Math. Soc., Vol. 156 (1971), \href{https://doi.org/10.1090/S0002-9947-1971-0275367-4}{73--83}.
		
		\bibitem{royden}
		Royden, H. L. and Fitzpatrick, P. M., \emph{Real Analysis}, Pearson, 4th Ed, \href{https://www.pearson.com/en-us/subject-catalog/p/real-analysis-classic-version/P200000006336/9780134689494}{2018}.
		
		\bibitem{Sierpinski}
		Sierpi\'nski, W., \emph{Sur les fonctions d’ensemble additives et continues.}, Fund. Math., Vol. 3 (1922), \href{https://doi.org/10.4064/fm-3-1-240-246}{240--246}.
		
	\end{thebibliography}
\end{document}